\documentclass[final]{siamart1116}

\newcommand*\figpath{./Figs}

\usepackage{lipsum}
\usepackage{amsfonts}
\usepackage{amssymb}
\usepackage{graphicx}
\usepackage{epstopdf}
\usepackage{algorithmic}
\usepackage{amsopn}
\usepackage{dsfont, bm}
\usepackage{subfigure}
\usepackage{dutchcal}
\usepackage{amsbsy}

\ifpdf
  \DeclareGraphicsExtensions{.eps,.pdf,.png,.jpg}
\else
  \DeclareGraphicsExtensions{.eps}
\fi

 \usepackage{color}
 
\usepackage[normalem]{ulem}
 \newcommand {\R} {{\mathds R}}
 
 \newcommand {\C} {{\mathds C}}

 \newcommand{\ie}{{\em i.\thinspace{}e. }}
 \newcommand{\eg}{{\em e.\thinspace{}g. }}

 \newcommand{\TF} {{\hbox{\tiny TF}}}
 \def\rtf{\rho_{\hbox{\tiny TF}}}
 \def\vec#1{\mathbf{#1}}

 \newcommand{\agrad}{\bA^t \cdot\bnabla}
 \def\rtf{\rho_{\hbox{\tiny TF}}}

 \def\pscal#1#2{{\langle #1, #2 \rangle}}

 \def\grdGn{G_n}

 \def\Real{\Re}
  
\def\bA{{\bf A}}

\def\x{{\bf x}}
\def\u{{\bf u}}
\def\d{{\bf d}}
\def\g{{\bf g}}

\def\RR{{\mathbb R}}
\def\CC{{\mathbb C}}

\def\M{{\mathcal M}}
\def\X{{X}}

\def\D{{\mathcal D}}
\def\T{{\mathcal T}}
\def\TuM{{\cal T}_u {\cal M}}
\def\TunM{{\cal T}_{u_n} {\cal M}}
\def\PG{P_{u_n,H_A}G}
\def\PGprev{P_{u_{n-1},H_A}G}
\def\Ru{{{\cal R}}_u}
\def\Run{{{\cal R}}_{u_n}}

\def\bnabla{\boldsymbol{\nabla}}

\def\argmin{\mathop{argmin}}

\newcommand{\TheTitle}{Computation of Ground States of the
  Gross-Pitaevskii Functional via Riemannian Optimization}
\newcommand{\TheAuthors}{I.~Danaila, B.~Protas}

\headers{Computation of Ground States via Riemannian
    Optimization}{\TheAuthors}

  \title{{Computation of Ground States of the Gross-Pitaevskii
      Functional via Riemannian Optimization} 
\thanks{Submitted to the      editors DATE. \funding{I.~D.~acknowledges the support through
        the ANR (France) grant ANR-12-MONU-0007-01 BECASIM (call
        ``Mod{\'e}les Num{\'e}riques'').  B.~P.~acknowledges the
        support through an NSERC (Canada) Discovery Grant.
        Computational resources were provided by CRIANN (Centre
        R{\'e}gional Informatique et d'Applications Num{\'e}riques de
        Normandie, France) under the project 2015001. The authors
        acknowledge the generous hospitality of the Fields Institute
        in Toronto during the Thematic Program on Multiscale
        Scientific Computing (January--April, 2016).  }}}

\author{
  Ionut Danaila\thanks{Laboratoire de Math{\'e}matiques Raphael Salem,
Universit{\'e} de Rouen Normandie, Rouen, France
    (\email{ionut.danaila@univ-rouen.fr}, \url{http://lmrs.univ-rouen.fr/Persopage/Danaila/}).}
  \and
  Bartosz Protas\thanks{Department of Mathematics \& Statistics, McMaster University, Hamilton, ON, Canada (\email{bprotas@mcmaster.ca}, \url{http://www.math.mcmaster.ca/bprotas}).}
}

\ifpdf
\hypersetup{
  pdftitle={\TheTitle},
  pdfauthor={\TheAuthors}
}
\fi

\begin{document}


\maketitle

\begin{abstract}
  In this paper we combine concepts from Riemannian Optimization
  [P.-A.~Absil, R.~Mahony, R.~Sepulchre, Optimization algorithms on
  matrix manifolds, Princeton University Press, 2008] and the theory
  of Sobolev gradients [J.~W.~Neuberger, Sobolev gradients, Springer
  2010] to derive a new conjugate gradient method for direct
  minimization of the Gross-Pitaevskii energy functional with
  rotation. The conservation of the number of particles in the system
  constrains the minimizers to lie on a manifold corresponding to the
  unit $L^2$ norm. The idea developed here is to transform the
  original constrained optimization problem to an unconstrained
  problem on this (spherical) Riemannian manifold, so that {fast}
  minimization algorithms can be applied {as alternatives to more
    standard constrained formulations}.  First, we obtain Sobolev
  gradients using an equivalent definition of an $H^1$ inner product
  which takes into account rotation.  Then, the Riemannian gradient
  (RG) steepest descent method is derived based on projected gradients
  and retraction of an intermediate solution back to the constraint
  manifold. Finally, we use the concept of the Riemannian vector
  transport to propose a Riemannian conjugate gradient (RCG) method
  for this problem. It is derived at the continuous level based on the
  ``optimize-then-discretize'' paradigm instead of the usual
  ``discretize-then-optimize'' approach, as this ensures robustness of
  the method when adaptive mesh refinement is performed in
  computations.  {We evaluate various design choices inherent in
    the formulation of the method and conclude with recommendations
    concerning selection of the best options.}  Numerical tests
  carried out in the finite-element setting based on Lagrangian
  piecewise quadratic space discretization demonstrate that the
  proposed RCG method outperforms the simple gradient descent (RG)
  method in terms of rate of convergence. {While on simple
    problems a Newton-type method implemented in the {\tt Ipopt}
    library exhibits a faster convergence than the (RCG) approach, the
    two methods perform similarly on more complex problems requiring
    the use of mesh adaptation. At the same time the (RCG) approach
    has far fewer tunable parameters. Finally,} the RCG method is
  extensively tested by computing complicated vortex configurations in
  rotating Bose-Einstein condensates, a task made challenging by large
  values of the non-linear interaction constant and the rotation rate
  {as well as by strongly anisotropic trapping potentials}.
\end{abstract}



\tableofcontents

\section{Introduction}
\label{sec:intro}

The rotating Bose-Einstein condensate (BEC) represents a
highly-controllable quantum system offering an ideal framework to
study quantized vortices at a macroscopic level. A rich variety of
vortex states, from a single vortex line to dense Abrikosov vortex
lattices and giant vortices, were experimentally observed and
extensively studied in the last two decades (\eg
\cite{BEC-physQT-2006-tsubota}).  A standard mathematical approach to
describe equilibrium configurations with quantized vortices in
rotating BEC is the minimization of the Gross-Pitaevskii (GP) energy
functional with rotation \cite{BEC-book-2003-pita,BEC-phys-2006-lieb}.
In addition to the global minima, the so-called ``ground states'',
local minima are also of interest as they represent excited, or
meta-stable, states which are more likely to be observed in
experiments \cite{BEC-physV-2004-bretin}.  The minimizers have the
form of complex-valued wavefunction fields dependent on the space
variable, resulting in an infinite-dimensional minimization problem.
The complexity of the minimization problem is further compounded by a
constraint imposed on the $L^2$ norm of the minimizers which reflects
the conservation of the number of atoms in the condensate. For the
mathematical properties of the GP energy with rotation and the
corresponding ground states we refer the reader to
\cite{BEC-book-2006-Aftalion,BEC-math-2001-Aftalion-Du,BEC-phys-2006-lieb,BEC-review-2013-Bao-KRM}.

In this paper we address the problem of direct minimization of the GP
energy functional with rotation when large nonlinear interaction
constants and high rotation frequencies are considered.  A number of
approaches to direct minimization of the GP energy have been proposed
based on various standard and emerging mathematical methods for
optimization problems in finite dimension: Optimal Damping Algorithm
\cite{BEC-CPCs-2007-dion-cances,BEC-CPCs-2014-Hohenester}, Newton-like
method based on Sequential Quadratic Programming (SQP)
\cite{BEC-CPCs-2013-Caliari}, Interior Point Method (Ipopt)
\cite{dan-2016-CPC}, Inertial Proximal Algorithm (iPiano)
\cite{BEC-numm-2016-Besse-opt} and regularized Newton method with
trust region \cite{BEC-baow-2016-Newton}.  Alternative approaches
which do not involve direct minimization of energy rely instead on the
solution of the corresponding Euler-Lagrange system which has the form
of a nonlinear eigenvalue problem. In the latter context, a wide
variety of classical integration and iterative techniques have been
employed such as Newton's method \cite{BEC-baow-2003-gs}, Runge-Kutta
\cite{BEC-CPCs-2013-Caplan} and continuation methods
\cite{BEC-numm-2010-Chang}, etc. {We also mention the
  ``deflated'' Newton's method, recently proposed in \cite{ckf18a},
  which represents a systematic approach to determine several distinct
  solution branches.}

Another class of approaches, pioneered in \cite{BEC-baow-2004-Du},
relies on a normalized gradient flow for the GP functional and became
popular due to their efficiency and ease of implementation (see also
the review papers
\cite{BEC-review-2006-bao,BEC-review-2013-Bao-KRM,BEC-review-2013-antoine-besse-bao,BEC-review-2014-Bao-ICM}).
These methods consist in first solving the gradient flow equation for
the minimization of an unconstrained energy followed by a
normalization of this ``predictor'' solution to bring it back to the
constraint manifold.  Solution of the gradient flow equation can be
viewed as a pseudo-time (or imaginary time) integration of the
time-independent GP equation.  Discretization of the gradient-flow
equations using a (natural) steepest descent method would result in a
very inefficient explicit Euler scheme for the (imaginary-) time
integration.  For this reason in \cite{BEC-baow-2004-Du} the
gradient-flow equation was solved using a semi-implicit backward Euler
scheme which proved even more efficient than the classical
Crank-Nicolson scheme.  The convergence of the original scheme
suggested in \cite{BEC-baow-2004-Du} was recently improved in
\cite{BEC-numm-2014-antoine-duboscq-JCP,BEC-CPCs-2016-antoine-levitt}
by using different discrete preconditioners.  It is interesting to
note that the gradient-flow equation for the GP functional has
structure similar to the complex-valued heat equation which makes it
amenable to solution with different classical time-integration schemes
such as Runge-Kutta-Fehlberg \cite{BEC-physV-2001-perez-gcM}, backward
Euler
\cite{BEC-baow-2004-Du,BEC-math-2001-Aftalion-Du,BEC-baow-2008-sym},
second-order Strang time-splitting
\cite{BEC-baow-2004-Du,BEC-math-2001-Aftalion-Du}, combined
Runge-Kutta-Crank-Nicolson scheme
\cite{dan-2003-aft,dan-2004-aft,dan-2005}, etc.

As regards the development of numerical methods, there are two main
paradigms, namely, ``optimize-then-discretize'' and
``discretize-then-optimize'', depending on whether the gradient
expressions are derived at the continuous or discrete level. In the
first case, the Sobolev gradient approach
\cite{BEC-book-2010-neuberger} represents the gradient-flow method
formulated with respect to a judiciously selected inner product in a
Hilbert space, rather than the classical $L^2$ inner product.  The
required gradients are obtained via the Riesz representation theorem.
When discretized, the Sobolev gradient approach can be also
interpreted as suitable preconditioning applied to the $L^2$ gradient
\cite{BEC-Sobolev-2002-farago}. However, the key advantage of working
with the ``optimize-then-discretize'' formulation is that the form of
this preconditioning is dictated by the functional (Sobolev) setting
of the problem and thus avoids the technically complicated search for
a good (and mesh-independent) discrete preconditioner. Sobolev
gradient methods were successfully applied to minimize the GP energy
in the presence of rotation in
\cite{BEC-physV-2001-perez-gcM,dan-2010-SISC}.

The purpose of this contribution is to develop and validate an
efficient computational approach to minimization of the GP energy by
combining the Sobolev gradient method with concepts of Riemannian
optimization \cite{ams08,s94b}. This allows us to transform the
original constrained optimization problem to an unconstrained problem
on a Riemannian manifold with a very simple structure which is
amenable to solution using the conjugate gradient approach. {We
  remark that while the ``Riemannian structure'' of an optimization
  problem may be exploited at various levels, in the present study we
  will focus solely on the basic concepts of ``retraction'' and
  ``vector transport'' describing how information travels along a
  manifold, and will not, in particular, consider endowing the
  constraint manifold with a Riemannian metric. In other words, we
  will assume that the constraint manifold is equipped with the metric
  induced by the embedding space.}  We begin by formulating a
Riemannian version of the Sobolev gradient approach in which the
retraction operation ensures that the norm constraint is satisfied at
all discrete times.  Then, convergence is accelerated using a
Riemannian version of the conjugate gradient method which relies on
the notion of the vector transport applied to the gradient and the
descent direction. Such approaches are already well established in the
context of problems formulated in finite dimensions \cite{manopt}, but
have received only {limited attention} in the context of problems
in infinite dimension. {Convergence of the Riemannian versions of
  the BFGS quasi-Newton approach and of the Fletcher-Reeves conjugate
  gradients method was established in \cite{rw12}, whereas some
  applications were considered in
  \cite{BEC-math-1997-Alouges,p04,aa09,mn13}}. To the best of our
knowledge, they represent a new direction as regards minimization of
the GP energy. {In our study, we carefully evaluate various
  design choices inherent in the formulation of the method and come up
  with recommendations concerning selection of the best options.
  Then,} we demonstrate that in combination with a flexible
finite-element discretization involving adaptive grid refinement
\cite{dan-2010-JCP,dan-2016-CPC}, the proposed approach outperforms
{a number of first-order techniques and performs on a par with a
  Newton-type method implemented in the {\tt Ipopt} library
  \cite{Ipopt-line-search-2005}}.

The structure of the paper is as follows: in the next section we state
the mathematical model describing minimization of the GP energy; in \S
\ref{sec:G} we recall the Sobolev gradient method with its projected
gradient variant, whereas the Riemannian gradient and conjugate
gradient methods are introduced in \S {\ref{sec:riemann}};
numerical discretization based on adaptive finite elements and its
software implementation are discussed in \S \ref{sec:numer};
{design choices to be explored in the formulation of the method
  are identified in \S \ref{sec:design};} in \S \ref{sec:res-manuf} we
use the method of the manufactured solutions to estimate the
{speeds} of convergence of the different approaches; in \S
\ref{sec:results} we compute a number of challenging BEC
configurations with vortices in various arrangements; discussions and
conclusions are deferred to \S \ref{sec:final}.

\section{Mathematical Model}
\label{sec:model}

The energy of a rotating homogeneous BEC at zero temperature is given
in terms of the Gross-Pitaevskii (GP) energy functional
\cite{BEC-book-2003-pita,BEC-phys-2006-lieb}. After applying standard
scaling and dimension reduction
\cite{BEC-book-2003-pita,BEC-review-2013-Bao-KRM}, its non-dimensional
form defined here on a two-dimensional domain $\D \subset \RR^2$
becomes
\begin{equation}
\label{eq-scal-energ}
\displaystyle {E}(u) = \int_{\D} \left[ \frac{1}{2} |\bnabla u|^2 + C_\text{trap}\, |u|^2 + \frac{1}{2} C_g |u|^4 - i C_\Omega\, u^* \agrad u \right] \, d{\vec{x}},
\end{equation}
where $u \; :\ \; \D \rightarrow \CC$ is a complex-valued wavefunction
and $u^*$ its complex conjugate, ${\bf A}^t = (y , -x)$,
$C_\text{trap} \; : \; \D \rightarrow \RR$ is the trapping potential.
$C_g$ and $C_{\Omega}$ are real constants characterizing the strength
of the nonlinear interactions and rotation frequency, respectively.
The wavefunction $u$ vanishes outside the trap and is therefore
assumed to satisfy the homogeneous Dirichlet boundary conditions $u =
0$ {on} $\partial\D$.  The conservation of the number of atoms in the
condensate is expressed as
\begin{equation}
\label{eq-scal-cons}
\| u \|_{2} := \| u \|_{L^2({\cal D},\CC)}={\sqrt{\int_{\D} \left|u({\vec x})\right|^2\, d{\vec x}}} = 1
\end{equation}
and serves as a constraint on $u$.  For the energy functional
\eqref{eq-scal-energ} to be well-defined, the wavefunction $u$ must
belong to the Sobolev space $H_0^1(\D,\CC)$ of functions with
square-integrable gradients \cite{af05} and vanishing traces on the
boundary (precise definitions of the norms in this function space will
be provided in \S \ref{sec:G} below). The constraint
\eqref{eq-scal-cons} may now be interpreted as defining a manifold
$\M$ in the solution space, \ie
\begin{equation}
{\cal M} := \left\{ u \in H^1_0({\cal D},\CC) \; : \; \| u \|_{2} = 1 \right\}.
\label{eq:M}
\end{equation}
We assume the trapping potential to have the following general form
allowing us to represent different trapping potentials used in
experiments
\begin{equation}
\label{eq-scal-trap-V}
C_\text{trap}(x,{y})=
\frac 12 \left(a_x x^2 + a_y y^2 + a_4 r^4\right), \quad r^2=(x^2+y^2),
\end{equation}
for some $a_x,a_y,a_4 \in \RR$. Along with the energy
\eqref{eq-scal-energ}, another important integral quantity describing
the rotating BEC is the total angular momentum
\begin{equation}
\label{eq-scal-Lzint}
L := {L}_z = i  \int_{\D} u^* \agrad u \, d{\vec x}
\end{equation}
{which under the assumed homogeneous Dirichlet boundary
  conditions is {real-valued}.}

Global minimizers of the energy functional \eqref{eq-scal-energ}
defined through
\begin{equation}
u_g = \arg\min_{u \in {\cal M}} E(u), \quad E(u_g) < \infty,
\label{eq-Emin}
\end{equation}
are called ground states.  Local minimizers, with energy larger than
that of the ground state, are referred to as excited or meta-stable
states.  The Euler-Lagrange system corresponding to the minimization
problem \eqref{eq-Emin} is derived using standard techniques and leads
to the stationary Gross-Pitaevskii equation
\begin{subequations}
\label{eq-scal-GP-stat}
\begin{alignat}{3}
-\frac{1}{2} \nabla^2 u + C_\text{trap} u + C_g |u|^2 u - i C_\Omega \agrad u &=&\,\, {\mu}\, u&
\quad & \text{in} \ \D, \label{eq-scal-GP-stat_a} \\
u & =& 0 & & \text{on} \ \partial\D, \label{eq-scal-GP-stat_b} \\
{\| u \|_{2}} & =& 1&. & \label{eq-scal-GP-stat_c}
\end{alignat}
\end{subequations}
The ground state and excited states are therefore eigenfunctions of
the nonlinear eigenvalue problem \eqref{eq-scal-GP-stat}.

\section{Gradient Flows and Steepest Descent Sobolev Gradient Methods}
\label{sec:G}

Numerical techniques for the solution of optimization problem
\eqref{eq-Emin} can be derived from {a form of the gradient-flow
  equation which for practical reasons we state here in terms of the
  gradient of the energy functional \eqref{eq-scal-energ} rather than
  the gradient of the corresponding Lagrangian functional}
\begin{equation}
\begin{aligned}
\frac{d u}{d t} & = - \nabla_{\X} E(u), \qquad t>0, \\
u(0) & = u_0,
\end{aligned}
\label{eq:Gflow}
\end{equation}
where $u_0 \in \M$ is an initial guess and $\nabla_{\X} E(u)$
represents the gradient of the GP energy functional
\eqref{eq-scal-energ} at $u$, computed with respect to the topology of
the Hilbert space $\X$ (to be made specific below). The gradient flow
needs to be additionally constrained to ensure that $u(t) \in \M$ for
$t>0$. {This approach is similar to the so-called {\em
    normalized gradient-flow method} \cite{BEC-baow-2004-Du}, which
  first {evolves} \eqref{eq:Gflow} and then projects the
  {intermediate} solution back onto the manifold. It can be
  viewed as a splitting method for solving the {{\em continuous
      normalized gradient-flow equation} \cite{BEC-baow-2004-Du} which
    is the} constrained {version of problem \eqref{eq:Gflow} in
    which the gradient $\nabla_{\X} E(u)$ is replaced with the
    gradient of the corresponding Lagrangian.}}

As shown below, many different computational approaches can be
derived from \eqref{eq:Gflow} by making specific choices of (i) the
Hilbert space $\X$, (ii) discretization of the initial-value problem
\eqref{eq:Gflow} with respect to pseudo-time $t$ and (iii) how the
constraint $u(t) \in \M$ is imposed.

{As regards the expression of the gradient $\nabla_{\X} E(u)$, it
  can be derived {from the G\^ateaux differential of energy
    \eqref{eq-scal-energ} using} the Riesz representation theorem
  \cite{l69} {which depends on the choice of the inner-product
    space $X$.}  Since energy \eqref{eq-scal-energ} is a twice
  continuously differentiable function from $H_0^1(\D,\CC)$ to $\RR$,
  a natural choice for the inner product that will ensure the
  existence of a gradient is}
\begin{equation}
	\pscal{u}{v}_{H^1}  = \int_{\mathcal D}  {(u,v)  +  (\bnabla u,\bnabla v)}\, d\x, 
	\label{eq-ipH1} 
\end{equation}
where $( u, v ) = u v^*$ is the {complex ($\C$ or $\C^2$) inner
  product}. A new inner product equivalent (in the precise sense of
norm equivalence) to \eqref{eq-ipH1} was suggested in
\cite{dan-2010-SISC}
\begin{equation}
\pscal{u}{v}_{H_A} = \int_{\mathcal D}  {(u,v)  +  (\bnabla_A u,\bnabla_A v)}\, d\x, \quad \bnabla_A = \bnabla + i
C_\Omega \bA, 
\label{eq-ipHA}
\end{equation}
and will be adopted in our considerations below.  Its definition was
motivated by the following physically more revealing form of the
energy functional equivalent to \eqref{eq-scal-energ}
 \begin{equation}
\label{eq-scal-energ-Agrad}
\displaystyle {E}(u) = 
\int_{\cal D} \left[ \frac{1}{2} \left|\bnabla u + i C_\Omega\, \bA u\right|^2 + C^\text{eff}_\text{trap}\, |u|^2 + \frac{1}{2} C_g |u|^4 \right] \, d{\vec x},
\end{equation}
where the effective non-dimensional trapping potential
$C^\text{eff}_\text{trap}$ is obtained from the original potential by
subtracting a term representing the centrifugal force
\cite{BEC-phys-1999-Stringari-TF}
\begin{equation}
\label{eq-scal-Ctrap-eff}
C^\text{eff}_\text{trap} (x,y)=C_\text{trap}(x,y)-\frac{1}{2} C_\Omega^2 \, r^2.
\end{equation}
{We add that, since the solution space $H^1_0(\D, \CC)$ is a
  subspace of both $H^1(\D, \CC)$ and $H_A(\D, \CC)$, we will assume
  $H^1_0(\D, \CC)$ to be equipped with the inner product
  \eqref{eq-ipH1} or \eqref{eq-ipHA}, and the notation $X=H^1$ or
  $X=H_A$ will make it clear which one is used.}

{For each $u \in X$, one can find an element of $X$ denoted
  $\nabla_X E(u)$, such that the directional G\^ateaux derivative of
  the energy at $u$ in the direction $v$, which is a continuous linear
  functional from $\X$ to $\RR$, is expressed as
\begin{equation}
E'(u)v = \Real\left(\langle \nabla_X E(u), v \rangle_X\right), \quad
\forall v \in X,
\label{eq-Riesz}
\end{equation}
where $\Real(\cdot)$ denotes the real part of a complex number.  We
refer to such an element of $X$ as the gradient of $E$ at $u$.
Computing the G\^ateaux derivative of the energy functional
\eqref{eq-scal-energ}, we obtain
\begin{equation}
E'(u)v= 2\Real\left( \int_{\mathcal D} \left[\frac{1}{2} (\bnabla u, \bnabla v) + \left( C_\text{trap} u + C_g \lvert u \rvert ^2 u - iC_\Omega \bA^t\cdot \bnabla u, v \right) \right]\, d{\vec x} \right)
\label{eq-derivgp}
\end{equation} 
which, together with \eqref{eq-Riesz}, allows us to identify the gradient
$\nabla_X E(u)$. In particular, the $H_A$ gradient, hereafter denoted 
$G = \nabla_{H_A} E(u)$, will be computed by solving the elliptic
boundary-value problem resulting from \eqref{eq-derivgp},
\eqref{eq-Riesz} and \eqref{eq-ipHA}, which we state here in the
equivalent weak form
\begin{multline}
\label{eq-num-algoS1}
{\forall v \in {H_0^1(\D, \R)}},  \int_{\mathcal D} \left[\left( 1 + C_\Omega^2 (x^2 + y^2) \right) (G, v) + (\bnabla G, \bnabla v) - 2i C_\Omega (\bA^t \cdot\bnabla G, v)\right]\, d{\vec x}  \\
= 2 \int_{\mathcal D} \left[\frac{1}{2} (\bnabla u, \bnabla v) + \left( C_\text{trap} u + C_g \lvert u \rvert ^2 u - iC_\Omega \bA^t\cdot \bnabla u, v \right) \right]\, d{\vec x}.
\end{multline}
}

\subsection{Normalized gradient flow}
\label{sec:Ngradflow}

{We note from \eqref{eq-Riesz} and \eqref{eq-derivgp} that, in
  order for $E'(u)$ to be bounded in the $L^2(\D,\CC)$ norm, stronger
  regularity assumptions must be imposed on the solution {$u$,
    namely $u \in H^2(\D,\CC)$}. In this case, from \eqref{eq-derivgp}
  we obtain that
\begin{equation}
E'(u)v= 2\Real\left( \int_{\mathcal D} \left({-\frac{1}{2}}  \nabla^2 u + C_\text{trap} u + C_g \lvert u \rvert ^2 u - iC_\Omega \bA^t\cdot \bnabla u, v \right) \, d{\vec x} \right),
\label{eq-derivgpL2}
\end{equation} 
{which allows us to formally} derive a "$L^2$-gradient"
corresponding to the $L^2$ inner product
\begin{equation}
	\pscal{u}{v}_{L^2}  = \int_{\mathcal D}  (u,v)\, d\x.
	\label{eq-ipL2}
\end{equation}
{We add that this is the expression appearing on the left-hand
  side of the Euler-Lagrange equation \eqref{eq-scal-GP-stat_a}.}}
The gradient flow equation \eqref{eq:Gflow} with this $L^2$-gradient
was discretized in \cite{BEC-baow-2004-Du} using a semi-implicit
backward Euler (BE) method
\begin{equation}
\frac{{\tilde u}-u_n}{\delta t} = \frac{1}{2} \nabla^2 {\tilde u} - C_\text{trap} {\tilde u} - C_g |u_n|^2 {\tilde u} + i C_\Omega \agrad {\tilde u},
\label{eq-grad-flowBE}
\end{equation}
where $u_n = u(t_n)$ denotes the approximation obtained at the $n$th
discrete time level, $\tilde{u} = \tilde{u}(t_{n+1})$ is an
intermediate (predictor) field and $\delta t$ is a fixed
{(pseudo-)}time step. The approximation $u_{n+1}$ at the time
level $t_{n+1}$ has to satisfy the unit-norm constraint
\eqref{eq-scal-cons} and is therefore obtained by normalizing the
predictor solution as
\begin{equation}
u_{n+1}= \frac{{\tilde u}(t_{n+1})}{\|{\tilde u}(t_{n+1})\|_{2}}.
\label{eq-steep-normB}
\end{equation}
This approach is referred to as the {\em normalized gradient flow}
method (see also
\cite{BEC-review-2006-bao,BEC-review-2013-Bao-KRM,BEC-review-2013-antoine-besse-bao,BEC-review-2014-Bao-ICM}).
Different existing variants of this method use {various}
numerical approaches to integrate the gradient-flow equation
\eqref{eq:Gflow}, \eg Runge-Kutta methods
\cite{BEC-physV-2001-perez-gcM,dan-2003-aft,dan-2004-aft,dan-2005},
different implicit schemes
\cite{BEC-baow-2004-Du,BEC-math-2001-Aftalion-Du,BEC-baow-2008-sym}
and Strang-type time-splitting approaches
\cite{BEC-baow-2004-Du,BEC-math-2001-Aftalion-Du}. Even though some of
these schemes do possess the energy-decreasing property, they
typically do not preserve the gradient-flow structure at the discrete
level, {in the sense that the expression on the right-hand side
  (RHS) of the gradient-flow equation \eqref{eq:Gflow} discretized
  with respect to the pseudo-time $t$ is no longer in the form of a
  gradient of $E(u)$ (this is because, as a result of the hybrid
  explicit/implicit treatment, different terms in this expression may
  depend on the state $u$ approximated at different time levels, as
  happens in \eqref{eq-grad-flowBE})}. Therefore, such imaginary-time
methods can be regarded as solving the nonlinear eigenvalue problem
\eqref{eq-scal-GP-stat}, rather than directly minimizing the GP
energy. Another potential drawback of such approaches is that
solutions of \eqref{eq-scal-GP-stat} are in general critical points of
the energy which are not necessarily minima.

\subsection{Steepest Descent Sobolev Gradient Methods}

Hereafter we focus on techniques which do preserve the gradient-flow
structure of \eqref{eq:Gflow} on the discrete level while explicitly
accounting for the presence of the unit-norm constraint
\eqref{eq-scal-cons}. As a starting point, we will thus consider an
explicit discretization of \eqref{eq:Gflow} in the following generic
form
\begin{equation}
\label{eq-num-descent}
u_{n+1} = u_{n} - \tau_n \, \grdGn, \quad n=0,1,\dots, 
\end{equation}
where $\tau_n$ is a suitable descent step size, whereas $\grdGn =
G(u_n)= \nabla_\X E(u_n)$ is a Sobolev gradient defined for
{$\X=H_A$} or $H^1$.  Below we discuss two ways in which the
information about the constraint $u \in \M$ can be incorporated in the
gradient method.

\subsubsection{Projected Sobolev Gradient Method}

\begin{figure}[!h]
	\begin{center}
		\includegraphics[width=0.8\textwidth]{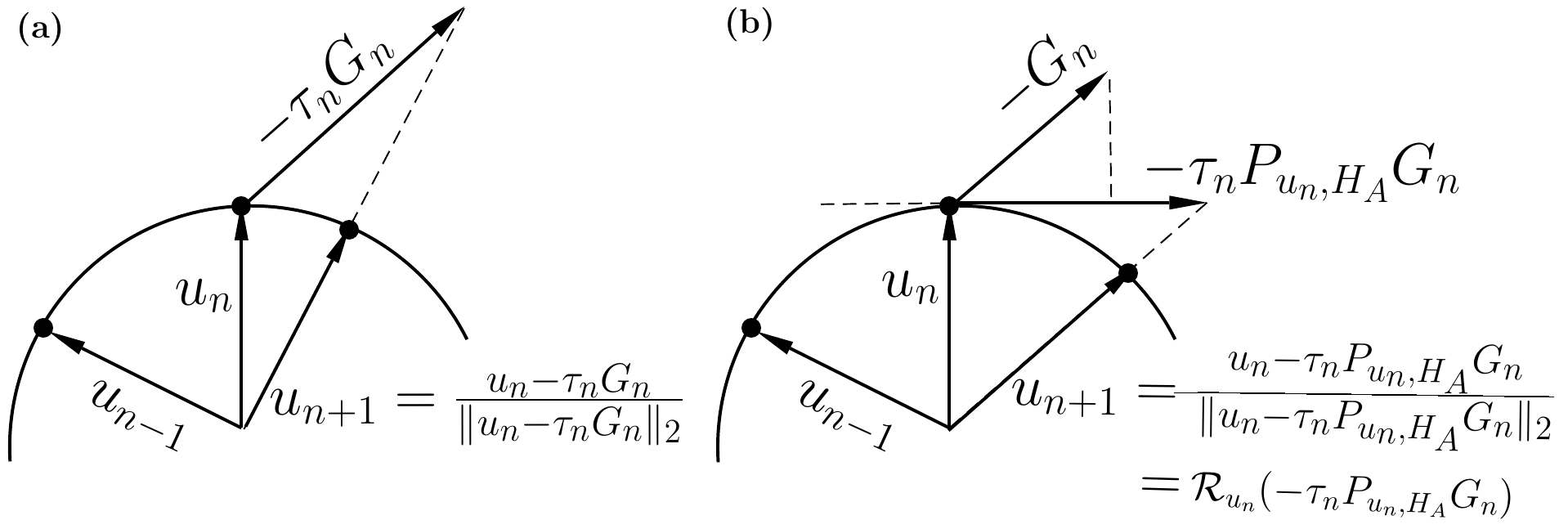}
	\end{center}
	\caption{Schematic illustration of the principle of the
          steepest descent method on the spherical manifold $\cal M$
          for (a) the simple (``unprojected'') gradient method and (b)
          the projected gradient (PG) {and the Riemannian
            gradient (RG) methods}.  In case (b), since the
          {projected gradient $P_{u_n, H_A} G_n$} belongs to the
          subspace $\TunM$ tangent to the manifold at $u_n$,
          normalization is equivalent to Riemannian retraction
          \eqref{eq:R}.}
	\label{fig:grad_proj}
\end{figure}

By considering the following identity derived from
\eqref{eq-num-descent}
\begin{equation}
\| u_{n+1} \|_{2}^2 = \| u_{n} - \tau_n \, \grdGn \|_{2}^2 =  \| u_{n}\|_{2}^2 - 2 \tau_n \Real\langle u_n, \grdGn\rangle_{L^2} + \tau_n^2  \| \grdGn \|_{2}^2,
\label{eq-norm-error}
\end{equation}
we note that using an unconstrained gradient $\grdGn$ leads to an
$\mathcal{O}(\tau_n)$ error in the satisfaction of the constraint
\eqref{eq-scal-cons} at each iteration.  Normalization of the solution
is then necessary to bring it back onto the manifold $\cal M$ (see
Figure \ref{fig:grad_proj}a).

This error can be reduced to second order $\mathcal{O}(\tau_n^2)$ by
requiring that $\langle u_n, \grdGn\rangle_{L^2}=0$, which is achieved
by projecting the gradient $\grdGn$ onto the subspace
\begin{equation}
\label{eq-TunM}
\TunM = \left\{ v \in H_0^1(\D,\CC) \; : \; \langle u_n, v \rangle_{L^2} = 0 \right\} 
\end{equation}
tangent to the constraint manifold $\cal M$ at $u_n$. As was shown in
\cite{dan-2010-SISC,BEC-numm-2009-kazemi}, the associated projection
operator $P_{u_n,\X}$ can be expressed in the following general form
\begin{equation}
  \label{eq-num-algoS3}
  P_{u_n,\X}G_n = G_n -  \lambda \, v_{\X},\qquad 
\lambda = \frac{ \Real \langle u_n , G_n \rangle_{L^2}}{ \Real \langle u_n , v_{\X} \rangle_{L^2}},
\end{equation}
where 
$v_{\X}$ is a solution of the variational problem
\begin{equation}
  \label{eq-num-algoS3b}
  \langle v_{\X} , v \rangle_{\X} = \langle u_n , v \rangle_{L^2}, \quad \forall v \in {\X}.
\end{equation}
We note that if $X=L^2$, $v_X=u$ in (\ref{eq-num-algoS3b}) and we
recover the well-known explicit expression of the projected
$L^2$-gradient (\eg \cite{BEC-math-1997-Alouges}).

Hereafter we will set $\X = H_A$ and denote $G_n = \nabla_{H_A}
E(u_n)$. Replacing $G_n$ with $P_{u_n, H_A}G_n$ in
\eqref{eq-num-descent}, we obtain the {\em projected gradient (PG)
  method} suggested in \cite{dan-2010-SISC}
\begin{equation}
\text{(PG)} \qquad 
u_{n+1} = u_n -  \tau_n \  P_{u_n,H_A}G_n, \quad n=0,1,\dots,
\label{eq:PG}
\end{equation}
While in \cite{dan-2010-SISC}
a fixed step size $\tau_n = \tau > 0$ was used, here  we use
an optimal step size found through the solution of a line-minimization
problem
\begin{equation}
  \label{eq-taun}
  \tau_n = \argmin_{\tau >0} \, E \left(u_n - \tau \, P_{u_n,H_A}G_n\right).
\end{equation}
An explicit expression for the optimal {descent} step was derived
in \cite{dan-2016-CPC} based on a particular form of the GP energy.
In this study, we prefer to solve problem \eqref{eq-taun} with a
general line-minimization approach such as Brent's algorithm
\cite{pftv86,nw00} as it has the advantage of being easily adapted to
the Riemannian gradient methods presented in the next section.  To
mitigate the $\mathcal{O}(\tau_n^2)$ drift away from the constraint
manifold $\M$ allowed by the (PG) iterations, normalization analogous
to \eqref{eq-steep-normB} may be applied to the iterates $u_n$ after a
certain number of steps. The idea of the (PG) approach is illustrated
schematically in Figure \ref{fig:grad_proj}b.

\section{Riemannian Optimization}
\label{sec:riemann}

In this section we discuss some basic concepts relevant to
optimization on manifolds, known as Riemannian optimization
\cite{ams08,s94b}. 
In contrast to the perspective developed in the previous section, here
we pursue a different, ``intrinsic'' approach where optimization is
performed directly on the manifold $\cal
M$. 
The main advantage of such a formulation is that it allows one to
treat \eqref{eq-Emin} as an {\em unconstrained} optimization problem
{creating an opportunity to apply a suitable modification of the
  conjugate gradient method as an alternative to more traditional
  constrained approaches such as, \eg {Sequential Quadratic
    Programming (SQP) for nonlinear optimization}.}

In addition to the definition of the projection on the tangent space
$\TuM$ already introduced above,
cf.~\eqref{eq-num-algoS3}--\eqref{eq-num-algoS3b}, we need to
introduce two more concepts, namely, the ``retraction'' (also referred
to as ``exponential mapping'') and the associated ``vector
transport''.  While in general these operators can have a rather
complicated form, in the present problem where the constraint manifold
$\cal M$ is given by \eqref{eq:M}, they {can be} reduced to
fairly simple expressions. We refer the reader to the monograph
\cite{ams08} for additional details concerning the
differential--geometry foundations of this approach.

\subsection{Riemannian Gradient Method}

Given a tangent vector $\xi \in \TuM$, where $u \in
{\cal M}$ is a state on the manifold, the {\em retraction} $\Ru \; :
\; \TuM \rightarrow {\cal M}$ is defined as the operator
\begin{equation}
\Ru(\xi) = \frac{u + \xi}{\| u + \xi \|_{2}},
\label{eq:R}
\end{equation}
where the norm used in the denominator is the same as the norm
defining the constraint manifold in \eqref{eq:M}. We note that for the
spherical manifold $\M$ the retraction operator is equivalent to
normalization \eqref{eq-steep-normB} already used in the previous
sections.

{\em The Riemannian gradient (RG) method} is then obtained applying
relation \eqref{eq:R} to the projected gradient $P_{u_n,H_A}G_n$,
cf.~\eqref{eq-num-algoS3}--\eqref{eq-num-algoS3b}, used in the (PG)
approach which yields
\begin{equation}
\text{(RG)} \qquad 
u_{n+1} = \Run\left( -\tau_n \PG_n \right), \qquad n=0,1,\dots.
\label{eq:RG}
\end{equation}
The step size $\tau_n$ is found optimally by solving a generalization
of the line-minimization problem \eqref{eq-taun} which uses retraction
\eqref{eq:R} to constrain the samples to manifold $\M$, \ie
\begin{equation}
\tau_n = \argmin_{\tau>0} E\left( \Run(-\tau \PG_n) \right).
\label{eq:tau}
\end{equation}
We refer to problem \eqref{eq:tau} as ``arc-minimization''. {It
  is solved} using a straightforward modification of Brent's algorithm
\cite{pftv86,nw00}. In addition to application of retraction
\eqref{eq:R} at every iteration in the latter case, the key difference
between the (PG) and (RG) approaches lies in how the optimal step size
$\tau_n$ is determined, cf.~\eqref{eq-taun} vs.~\eqref{eq:tau}. The
idea of the (RG) method is illustrated schematically in Figure
\ref{fig:grad_proj}b.

\subsection{Riemannian Conjugate Gradient Method}
\label{sec:RCG}

As a point of reference, we begin by recalling the nonlinear conjugate
gradients method in the Euclidean case. Given a function $f \; : \;
\R^N \rightarrow \R$, this approach finds its local minimum
$\overline{\u}$ as $\overline{\u} = \lim_{n \rightarrow \infty} \u_n$,
with the iterates $\u_n$ defined as follows \cite{nw00}
\begin{equation}
\u_{n+1} = \u_n + \tau_n \, \d_n, \qquad n=0,1,\dots,
\label{eq:CG}
\end{equation}
where $\u_0$ is the initial guess and the descent direction $\d_n$ is
constructed as
\begin{equation}
\begin{aligned}
\d_0 & = -\g_0, \\ 
\d_{n} & = -\g_n + \beta_n \, \d_{n-1}, \qquad n=1,2,\dots 
\end{aligned}
\label{eq:dn}
\end{equation}
in which $\g_n = \bnabla f(\u_n)$ and $\beta_n$ is a ``momentum'' term
chosen to enforce the conjugacy of the search directions $\d_k$,
$k=1,\dots,n$. When the objective function is convex-quadratic, \ie
$f(\u) = \u^T \bA \u$ for some positive-definite matrix $\bA \in
\R^{N\times N}$, approach \eqref{eq:CG}--\eqref{eq:dn} reduces to the
``linear'' conjugate gradient method in which $\tau_n$ and $\beta_n$
are given in terms of simple expressions involving $\bA$ \cite{nw00}.
In the non-quadratic setting, which is the case of problem
\eqref{eq-Emin}, the step size $\tau_n$ needs to be found via line
minimization as described by \eqref{eq:tau}, whereas the momentum term
is typically computed using one of the following expressions
\begin{subequations}
\label{eq:beta}
\begin{alignat}{2}
\beta_n &= \beta_n^{FR} := \frac{\langle \g_n, \g_n \rangle_{{\cal Y}}}{\langle \g_{n-1}, \g_{n-1} \rangle_{{\cal Y}}}& \qquad & \text{(Fletcher-Reeves)},
\label{eq:betaFR} \\
\beta_n &= \beta_n^{PR} := \frac{\langle \g_n, (\g_n - \g_{n-1})\rangle_{{\cal Y}}}{\langle \g_{n-1}, \g_{n-1} \rangle_{{\cal Y}}}& \qquad & \text{(Polak-Ribi{\`{e}}re)},
\label{eq:betaPR} 
\end{alignat}
\end{subequations}
where $\langle \cdot,\cdot \rangle_{{\cal Y}}$ is the inner product
defined with respect to the metric $\cal Y$ (in the simplest case when
$\u \in \R^N$, $\langle {\bf{ a}},{\bf{b}} \rangle_{{\cal Y}} = {\bf{
    a}}^T{\bf{b}}$ for ${\bf{ a}},{\bf{b}} \in \R^N$). The coefficient
$\beta_n$ {may be} periodically reset to zero which is
{known to improve} convergence for convex, non-quadratic problems
\cite{nw00}.  It is well known that for optimization problems which
are locally quadratic the conjugate gradient approach exhibits faster
(though still linear) convergence than the convergence characterizing
the simple gradient method, especially for poorly scaled
{problems} \cite{nw00}. Similar observations have been also
reported for Riemannian optimization problems \cite{ams08,s94b}.

We explain below how the conjugate gradients approach can be adapted
to the Riemannian case involving the energy functional
\eqref{eq-scal-energ} defined for infinite-dimensional state variables
$u \in H_0^1({\cal D})$.  There are two key issues which must be
addressed:
\begin{enumerate}
\item the two {terms} on the RHS of formula \eqref{eq:dn} belong
  to two different linear spaces which are the tangent spaces
  constructed at two consecutive iterations, \ie $\g_n \in {\cal
    T}_{\u_n} {\cal M}$ and $\d_{n-1} \in {\cal T}_{\u_{n-1}} {\cal
    M}$; as a result, they cannot be simply added; the same problem
  also concerns the inner-product expressions in the numerator of the
  Polak-Ribi{\`{e}}re momentum term \eqref{eq:betaPR},

\item while in the finite {dimensions} all norms are equivalent,
  this is no longer the case in the infinite-dimensional setting where
  the choice of the metric does play a significant role; in our
  approach, although the gradient-descent equations are discretized
  {in space} for the purpose of the numerical solution, their
  specific form is derived in the infinite-dimensional setting (in
  other words, we follow the ``optimize-then-discretize'' paradigm
  \cite{g03}); in addition to the momentum term \eqref{eq:beta}, the
  choice of the metric implied by the inner product also plays a role
  in the construction of the projection
  \eqref{eq-num-algoS3}--\eqref{eq-num-algoS3b} and the vector
  transport which will be defined below.
\end{enumerate}

The key concept required in order to address the first issue is the
{\em vector transport}
\begin{equation}
{{\cal TM} \oplus {\cal TM} \rightarrow {\cal TM} \; : \; 
(\eta,\xi) \longmapsto \T_{\eta}(\xi) \in {\cal TM},}
\label{eq:T}
\end{equation}
where ${\cal TM} = \cup_{x \in {\cal M}} {{\cal T}_{x} {\cal M}}$ is
the tangent bundle, describing how the vector field $\xi$ is
transported along the manifold $\cal M$ by the field $\eta$
\cite{ams08}. It therefore generalizes the concept of the parallel
translation to the motion on the manifold and is also closely related
to the ``affine connection'' which is one of the key
differential-geometric quantities characterizing a manifold.  The
vector transport thus provides a map between the tangent spaces
${{\cal T}_{u_{n-1}} {\cal M}}$ and $\TunM$ obtained at two
consecutive iterations, so that algebraic operations can be performed
on vectors belonging to these subspaces. 

In general, vector transport is not defined uniquely and in the
present case when the manifold is a sphere, the following two
definitions lead to expressions particularly simple from the
computational point of view. Let $u \in {\cal M}$ and $\eta_u, \xi_u
\in {\TuM}$; the transport of $\xi_u$ by $\eta_u$ can be expressed
either by
\begin{itemize}
	\item vector transport via differentiated retraction
	\begin{equation}
	\T_{\eta_u}(\xi_{u}) = \frac{d}{dt} {\cal R}_{u}({\eta_{u} + t \xi_{u}}) \big|_{t=0} = 
	\frac{1}{\| u + \eta_{u} \|}\left[ \xi_{u} - \frac{\langle u+\eta_{u}, \xi_{u} \rangle}{\| u + \eta_{u} \|^2} (u+\eta_{u}) \right],
	\label{eq:T1}
	\end{equation}
	\item or by vector transport {using Riemannian submanifold structure}
	\begin{equation}
	\T_{\eta_u}(\xi_{u}) = P_{{\cal R}_{u}(\eta_{u})} \xi_{u} = 
	\left[ \xi_{u} - \frac{\langle u+\eta_{u}, \xi_{u} \rangle}{\| u + \eta_{u} \|^2} (u+\eta_{u}) \right],
	\label{eq:T2}
	\end{equation}
\end{itemize}
where 
$P_u$ is the orthogonal projector on ${\TuM}$.  We note that these
formulas differ only by a scalar factor $\| u + \eta_u \|^{-1}$ and
{expression \eqref{eq:T2} can be interpreted as a Riemannian
  parallel transport. We further remark that} the vector transport
$\T_{\eta_u}(\xi_{u})$ is linear in the field $\xi_{u}$, but not in
$\eta_{u}$. {The reader is referred} to monograph \cite{ams08}
for details concerning the derivation of formulas
\eqref{eq:T1}--\eqref{eq:T2}.  The numerical results presented in
\S\ref{sec:res-manuf} and \S\ref{sec:results} are obtained using
{the vector-transport} {expression \eqref{eq:T1} or  \eqref{eq:T2}} with the $L^2$
inner product and norm. {This choice of the metric is dictated by
  the norm defining the constraint manifold,
  cf.~\eqref{eq-scal-cons}.}

\begin{figure}	
	\begin{center}
		\includegraphics[width=0.8\textwidth]{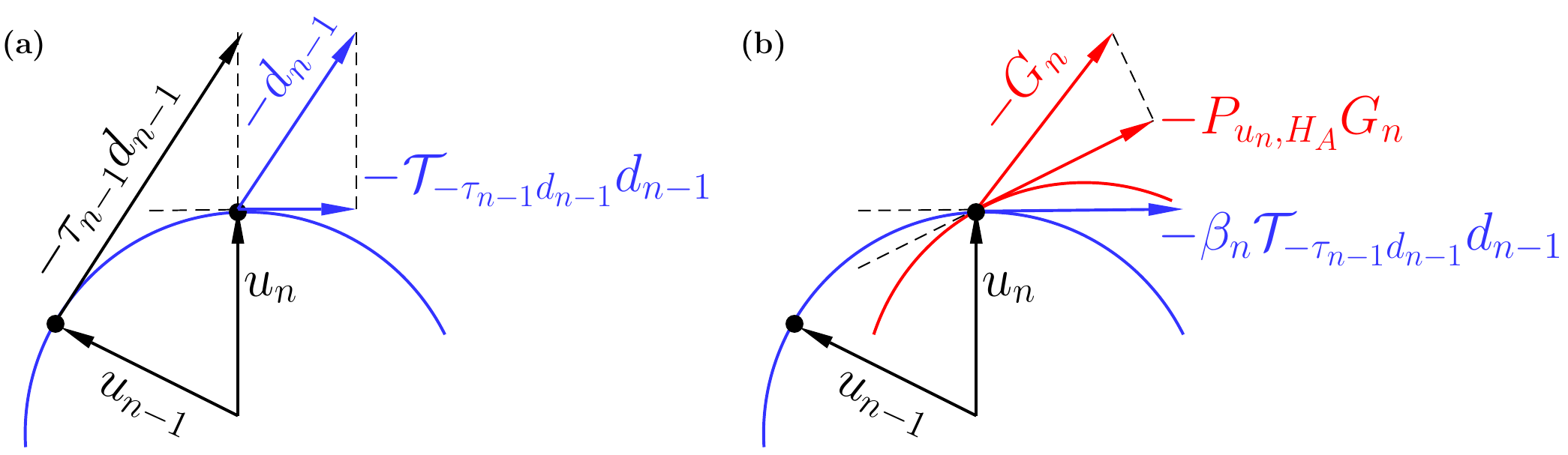}
	\end{center}
	\caption{Schematic illustration of the principle of Riemannian
          conjugate gradients (RCG) method on a spherical manifold.
          (a) Riemannian vector transport of the anterior conjugate
          direction $d_{n-1}$; the transport of the anterior gradient
          $G_{n-1}$ is performed in a similar way. (b) Projection of
          the new Sobolev gradient $G_n$ onto the tangent subspace
          $\TunM$ resulting in $\PG_n$.  The linear combination
          (\ref{eq:Rdn}) of $\PG_n$ and the transported anterior
          direction is computed in $\TunM$. }
	\label{fig:grad_conj}
\end{figure}

{Finally,} the conjugate gradient method
\eqref{eq:CG}--\eqref{eq:beta} can be rewritten in the Riemannian
infinite-dimensional setting as
\begin{equation}
\text{(RCG)} \qquad 
u_{n+1} = \Run\left(- \tau_n \, d_n \right), \qquad n=0,1,\dots,
\label{eq:RCG}
\end{equation}
where 
\begin{equation}
\begin{aligned}
d_0 & = - {P_{u_0,H_A}G_0}, \\
d_{n} & = -\PG_n  + \beta_n \, \T_{-\tau_{n-1} d_{n-1}} (d_{n-1}), \qquad n=1,2,\dots 
\end{aligned}
\label{eq:Rdn}
\end{equation}
with the Polak-Ribi{\`{e}}re momentum term modified as follows (the
corresponding term in the Fletcher-Reeves approach remains unchanged)
\begin{equation}
\beta_n = \beta_n^{PR} := \frac{\Big\langle \PG_n, \left(\PG_n - \T_{-\tau_{n-1} d_{n-1}} \PGprev_{n-1} \right)\Big\rangle_{H_A}}
{\Big\langle \PGprev_{n-1}, \PGprev_{n-1} \Big\rangle_{H_A}}.
\label{eq:RbetaPR} 
\end{equation}
The optimal descent step $\tau_n$ in \eqref{eq:RCG} is computed as in
\eqref{eq:tau} by solving the corresponding arc-minimization problem
\begin{equation}
\tau_n = \argmin_{\tau>0} E\left( \Run(-\tau d_n) \right)
\label{eq:tauRCG}
\end{equation}
using a generalization of Brent's method. We refer to approach
\eqref{eq:RCG}--\eqref{eq:tauRCG} as the {\em Riemannian Conjugate
  Gradients (RCG) method}, and its idea is schematically illustrated in
Figure \ref{fig:grad_conj}.

\section{Space Discretization}
\label{sec:numer}

We use a finite-element approximation constructed as follows.  Let
${T}_h$ be a family of triangulations of the domain $\cal D$
{parametrized by the mesh size $h>0$}.  We assume that ${T}_h$ is
a regular family {in the sense of Ciarlet \cite{ciarlet-1978}},
with {$h$} belonging to a generalized sequence converging to
zero. We denote by $P^l(T)$ the space of polynomial functions of
degree not exceeding $l\geq 1$ defined on triangles $T\in {T}_h$. We
also introduce the finite-element approximation spaces
\begin{eqnarray}
W_h^l &=& \left\{ w_h \in C^0({\bar {\cal D}_h}); \ w_h |_T \in P^l(T), \ \forall T\in {T}_h\right\},\\
V_h^l &=& \left\{ w_h \in W_h^l; \ w_h |_{\Gamma_h} = 0\right\}, \quad \text{where} \ \Gamma_h = \partial {\cal D}_h.
\end{eqnarray}
The finite dimensional space $V_h^l$ is a subspace of $H^1_0({\cal
  D})$ and therefore will be used to approximate the energy functional
\eqref{eq-scal-energ} and the different expressions representing
gradients and descent directions in algorithms (PG), (RG) and (RCG).
In the following we use {$P^4$ ($l=4$, piecewise quartic) finite
  elements to approximate the nonlinear terms in \eqref{eq-num-algoS1}
  and the $P^2$ representation for the remaining terms, for evaluation
  of the GP energy \eqref{eq-scal-energ} and also to represent the
  approximate solution $u_n$.} In addition, adaptive mesh refinement
suggested in \cite{dan-2010-JCP} and tested in
\cite{dan-2010-JCP,dan-2016-CPC} is used to adapt the grid during
iterations leading to a significant reduction of the computational
time.  The approach is implemented in {\tt FreeFEM++}
\cite{hecht-2012-JNM,freefem}, {where mesh adaptivity relies on
  metric control \cite{hecht-2000-ijnmf}. The main idea is to define a
  metric based on the Hessian and use a Delaunay procedure to build a
  new mesh such that all the edges are close to the unit length with
  respect to this new metric.  We use the adaptive meshing strategy
  suggested in \cite{dan-2010-JCP,dan-2016-CPC}. The relative change
  of the energy of the solution (cf.~\eqref{eq:stopping}) is used as
  an indicator to trigger mesh adaptation in which the metrics are
  computed simultaneously using the real and imaginary part of the
  solution.}  The implementation of the Riemannian retraction
\eqref{eq:R} and {vector transport \eqref{eq:T1} or
  \eqref{eq:T2}} is straightforward and was found to work very well
with arc-minimization \eqref{eq:tauRCG} and adaptive mesh refinement.

\section{{Design Choices Inherent in the (RCG) Method}}
\label{sec:design}

{As is evident from \S \ref{sec:RCG}, the (RCG) method offers a
  number of design choices which can be exploited to optimize its
  performance for a specific problem. One of the goals of this study is
  to evaluate these options in the context of minimization of the GP
  energy, and we will focus on the following choices most relevant for
  the Riemannian aspect of the proposed approach:
\begin{itemize}
\item[(i)] form of the momentum term $\beta_n$: Fletcher-Reeves
  \eqref{eq:betaFR} or Polak-Ribi{\`{e}}re \eqref{eq:betaPR}, with the
  corresponding variants of the (RCG) method referred to as (RCG)-(FR)
  and (RCG)-(PR), respectively;

\item[(ii)] form of the vector transport $\T_{\eta_u}(\xi_{u})$:
  defined via differentiated retraction \eqref{eq:T1} or using the
  Riemannian submanifold structure \eqref{eq:T2}, with the
  corresponding variants of the (RCG) method referred to as
  (RCG)-(VtDR) and (RCG)-(VtRS), respectively; in addition, we will
  also consider the classical conjugate gradients (CG) method without
  vector transport (\ie with $\T_{\eta_u}(\xi_{u}) = \xi_u$).
\end{itemize}
Combining these different choices yields already six distinct
variants, \ie (RCG)-(FR)-(VtDR), (RCG)-(FR)-(VtRS), (RCG)-(FR),
(RCG)-(PR)-(VtDR), (RCG)-(PR)-(VtRS) and (RCG)-(PR). Evidently, there
also exist other design choices which, however, will not be considered
here, because they are less relevant for the Riemannian aspect of the
problem and/or have already been considered elsewhere before. For
example, the choice of the metric $X$ defining the gradient in
\eqref{eq-Riesz} and the projection in
\eqref{eq-num-algoS3}--\eqref{eq-num-algoS3b} was extensively analyzed
in \cite{dan-2010-SISC,dan-2010-JCP,dan-2016-CPC}, and here we will
exclusively use $X = H_A$ which was found to be the best choice for
minimization of the GP energy in the presence of rotation.  For the
Riemannian operators of retraction and vector transport we use the
$L^2$ metric naturally induced by the spherical manifold defined in
\eqref{eq:M}.}

{In principle, we could also consider the frequency of
  retractions as yet another design parameter, however, since this
  operation has a negligible cost, it is performed at every iteration
  (which is also consistent with the need to reinterpolate the
  solution once a grid adaptation has been taken place). On the other
  hand, we will consider the effect of periodically resetting the
  momentum term $\beta_n$ to zero (cf.~\S \ref{sec:RCG}). In the
  sections to follow we will analyze these different design choices in
  order to identify the most robust (RCG) method capable of handling
  mesh adaptivity, which in earlier studies
  \cite{dan-2010-JCP,dan-2016-CPC} was shown to be indispensable for
  computational efficiency.}

\section{Convergence {Speeds} of Different Gradient Methods}
\label{sec:res-manuf}

\begin{figure}[!h]
	\begin{center}
		\includegraphics[width=0.6\textwidth]{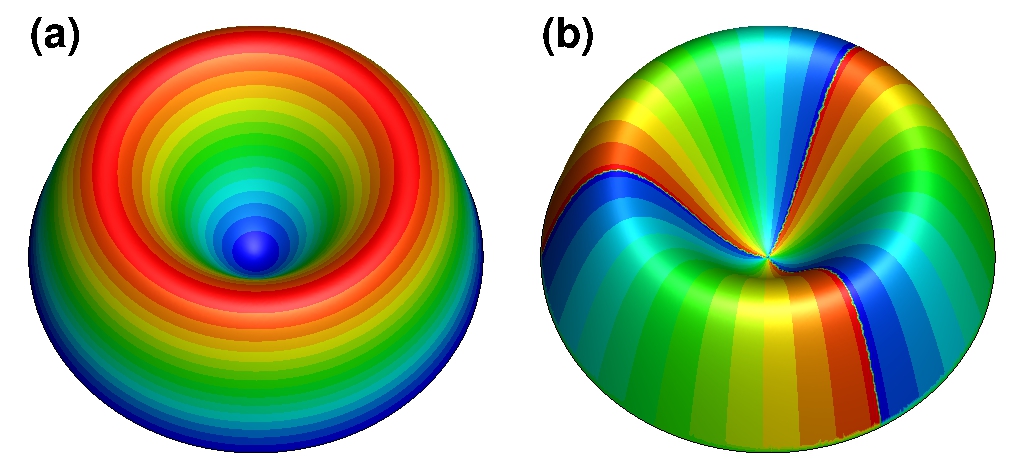}
	\end{center}
	\caption{Manufactured solution (\ref{eq-manuf-exact})
          visualized with a 3D-rendering of the modulus $|u_{ex}|$
          color-coded with (a) the modulus itself and (b) the phase of
          the solution for $m=3$.}
	\label{fig:manuf-exact}
\end{figure}

{We start by assessing the convergence speed of the minimization
  algorithms (PG), (RG), (CG) and {the different variants of
    the (RCG) approach} using {fixed meshes with different
    resolutions}. We use the method of manufactured solutions
  \cite{BEC-book-1998-Roache}, a general tool for verification of
  calculations which has the advantage of providing an exact solution
  to a modified problem, related to the {true} one.}  {The
  general idea is that,} by introducing an extra source term, the
original system of equations is modified to admit an exact solution
given by a convenient analytic expression. Even though in most cases
exact solutions constructed in this way are not physically realistic,
this approach allows one to rigorously verify computations.  Here we
manufacture such an exact solution in the form
\begin{equation}
u_{ex}(x,y) = U(r) \, \exp(i m \theta), \qquad U(r) 
= \frac{2 \sqrt {21}}{\sqrt {\pi}}\,{\frac {{r}^{2} \left( R-r \right) }{R^{4}}}, \quad m \in \mathbb{N},
\label{eq-manuf-exact}
\end{equation}
where $(r,\theta)$ are the cylindrical coordinates of the point
$(x,y)$ and $R$ is the radius of the circular domain $\cal D$. We note
that this solution satisfies constraint \eqref{eq-scal-cons} and
qualitatively resembles a giant vortex in the condensate (see Figure
\ref{fig:manuf-exact} and \S \ref{sec:results}).
It also satisfies an inhomogeneous form of the nonlinear problem
\eqref{eq-scal-GP-stat}, \ie
\begin{equation}
\label{eq-manuf-GP}
\begin{alignedat}{3}
-\frac{1}{2} \nabla^2 u_{ex} + C_\text{trap} u_{ex} + C_g |u_{ex}|^2 u_{ex} - i C_\Omega \agrad u_{ex} &=& f&
\quad & \text{in} \ \D, \\
u_{ex} & =& 0& & \text{on} \ \partial\D,
\end{alignedat}
\end{equation}
and is a critical point of  the modified energy functional
\begin{equation}
E(u,f) =
\int_{\D} \left[ \frac{1}{2} |\bnabla u|^2 + C_\text{trap}\, |u|^2 + 
\frac{1}{2} C_g |u|^4 - i C_\Omega\, u^* \agrad u - (f^* u + f u^*) \right] \, d{\vec{x}}.
\label{eq-manuf-Energy}
\end{equation}
For this energy functional, the $L^2$ gradient is expressed as
{discussed in \S \ref{sec:Ngradflow}}, but with a supplementary
term $-2f$ added.  Given the form \eqref{eq-manuf-exact} of the
manufactured solution and assuming a harmonic trapping potential
$C_{\text{trap}} =r^2/2$, from \eqref{eq-manuf-GP} we obtain
$f(r,\theta) =F(r) \, \exp(i m \theta)$, where $F(r)$ is a polynomial
of degree 9. From this and relations \eqref{eq-manuf-Energy} and
\eqref{eq-scal-Lzint} we can deduce exact expressions for the energy
$E_{ex} := E(u_{ex})$ and the angular momentum $L_{ex} :=
L(u_{ex})=m$.

The numerical tests are based on the manufactured solution
\eqref{eq-manuf-exact} corresponding to the following parameter
values: $C_{\text{trap}}=r^2/2$, $C_g=500$, $R=1$, $m=3$,
$C_\Omega=10$, where, to make the problem more challenging, large
values of the nonlinear interaction constant $C_g$ and rotation
frequency $C_\Omega$ are used (cf.~Figure \ref{fig:manuf-exact}).  The
step size $\tau_n$ is determined at each iteration by
{line-minimization \eqref{eq-taun} for the (PG) method and
  arc-minimization \eqref{eq:tau} or \eqref{eq:tauRCG} for the (RG),
  (CG) and (RCG) methods}.

In order to assess the mesh-independent effect of the Sobolev gradient
preconditioning, cf.~\S \ref{sec:numer}, we perform computations using
two grids: Mesh 1 consisting of 24,454 triangles with $h_{min}=0.0118$
and Mesh 2 consisting of 99,329 triangles with $h_{min}=0.0059$,
where $h_{min}$ is the smallest grid size. In the present case no mesh
refinement was performed during iterations. The initial guess is taken
as $u_0=0$, whereas iterations are declared converged once the
following condition based on the relative energy decrease is satisfied
\cite{dan-2010-JCP,dan-2016-CPC}
\begin{equation}
\varepsilon_E =|E_{n+1}-E_n|/E_n < \varepsilon_{st}=10^{-12}, \qquad \text{where} \ E_n := E(u_n).
\label{eq:stopping}
\end{equation}

{The performance of the approaches corresponding to the different
  design choices discussed in \S \ref{sec:design} is summarized in
  Table \ref{tab:manuf-perf} where all computations were performed on
  Mesh 2.  In this and in the following tables the CPU time reflects
  the computations performed on {a Linux workstation with two
    3.10GHz Intel Xeon E5-2687w CPUs.} 
\begin{table}[h]
	\begin{center}
		\begin{tabular}{|r||r|l|}\hline
			Method & iter & CPU \\ \hline \hline
			(RCG)-(PR)-(VtRS) & 37 &1270 \\ \hline
			(RCG)-(PR)-(VtDR) & 38 &1326 \\ \hline
			(CG)-(PR)             & 38 &1529 \\ \hline \hline
			(RCG)-(FR)-(VtRS)  & 54 &1852 \\ \hline
			(RCG)-(FR)-(VtDR)  & 49 &1668 \\ \hline
			(CG)-(FR)             & 31 &1297 \\ \hline \hline	
			(RG)                       & 180& 5274 	\\ \hline	
			(PG)                       &219& 3104 \\ \hline	
		\end{tabular}
	\end{center}
	\caption{{Test case based on the manufactured solution 
			\eqref{eq-manuf-exact}.  Performance of {the gradient 
				methods corresponding to the different design choices, 
				cf.~\S \ref{sec:design} measured} in terms of the number of iterations 
			(iter) and the computational time in seconds (CPU) 
			{required for} convergence. Note that the total computational time depends both on the number of iterations and the number of evaluations of the energy functional in the line- or arc-minimization procedures.}}
	\label{tab:manuf-perf}	
\end{table}

In the calculations reported in
  Table \ref{tab:manuf-perf} for the (CG) and (RCG) methods we did not
  reset the momentum term $\beta_n$ to zero.}  {Firstly, we note
  that all {variants of the} (RCG) and (CG) methods outperform
  the simple gradient methods (PG) and (RG). The reduced CPU time of
  the (PG) method comes from the fact that it implements a
  line-minimization {strategy} \eqref{eq-taun} based on
  analytical expressions derived from the particular form of the GP
  energy (see \cite{dan-2016-CPC} for details). This is not the case
  for the arc-minimization used by all Riemannian gradient methods,
  where {an adaptation of} Brent's algorithm is employed.
  Secondly, we observe that the (RCG) algorithm with the
  Polak-Ribi{\`{e}}re (PR) {momentum term is least sensitive to
    the form of the vector transport. On the other hand,} the
  Fletcher-Reeves version of the (RCG) method proves more sensitive to
  the form of the vector transport {and, somewhat surprisingly,
    the (CG)-(FR) approach (without vector transport) turns out to be
    the most efficient in terms of the number of iterations (although
    not in terms of the computational time). The method which
    converged in the shortest time was the (RCG)-(PR)-(VtRS) approach
    and it will be used for further tests in the remainder of this
    section; for brevity, we will refer to it simply as (RCG).}}

\begin{figure}[h]
  \begin{center}
    \includegraphics[width=0.75\textwidth]{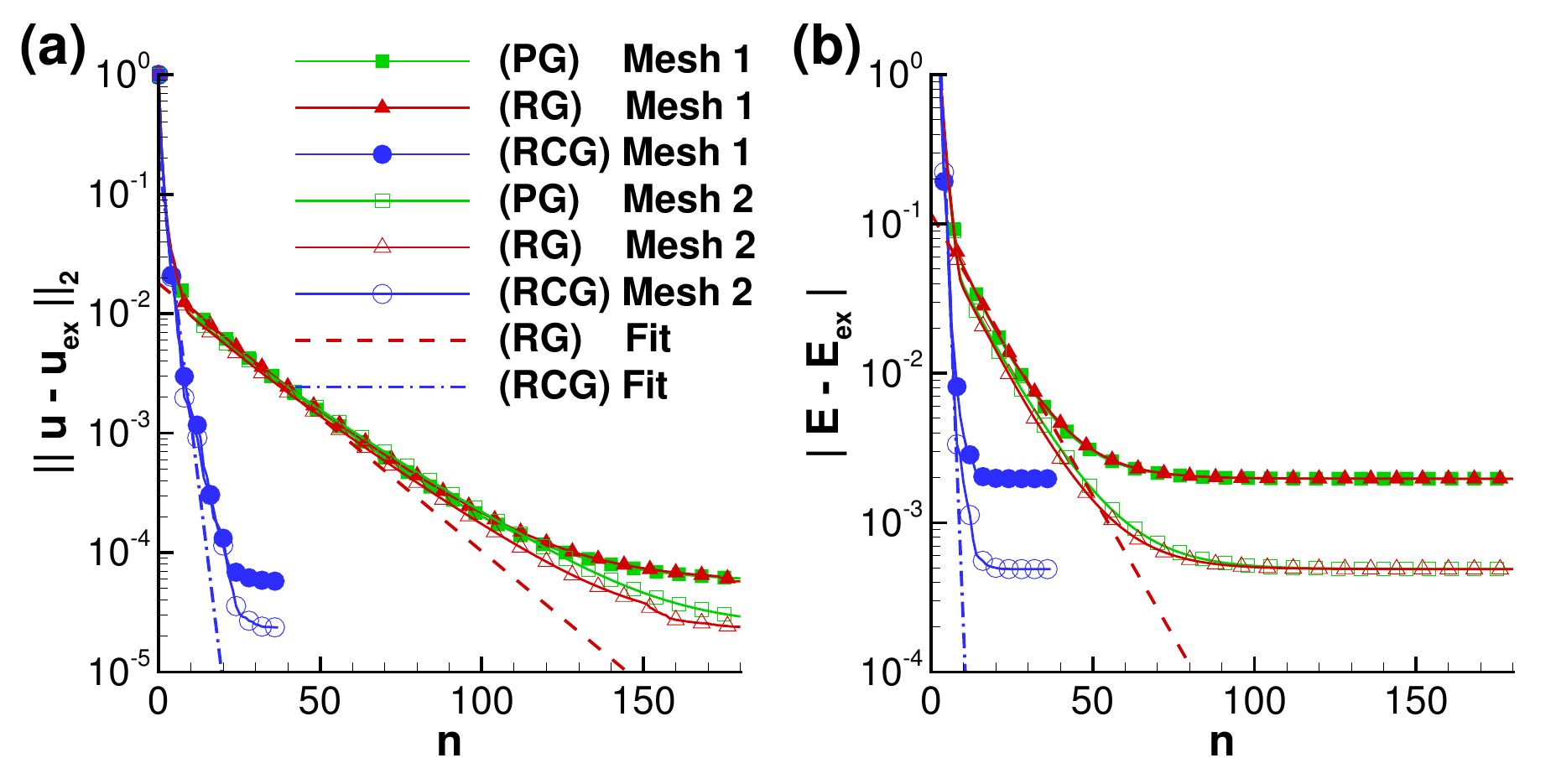}
  \end{center}
  \caption{Test case based on the manufactured solution
    \eqref{eq-manuf-exact}. Convergence of the (PG), (RG) and (RCG)
    methods: (a) $\|(u_{n}-u_{ex})\|_{2}$ and (b) $|E_{n}-E_{ex}|$ are
    shown as functions of iterations $n$ for different
    discretizations.  Dashed and dash-dotted lines indicate the
    least-squares fits \eqref{eq:manuf-exact-fit}.}
  \label{fig:manuf-exactE}
\end{figure}

To assess the {speed} of convergence of the (PG), (RG) and (RCG)
approaches, the quantities $\|u_n-u_{ex}\|_{2}$ and $|E_n-E_{ex}|$ are
shown as functions of iterations $n$ for the two spatial
discretizations in Figures \ref{fig:manuf-exactE}a and
\ref{fig:manuf-exactE}b.  In these figures we observe linear
convergence followed by a slower convergence at final iterations. The
change of the slope of error curves occurs at the level at which the
minimization errors $(u_n-u_{ex})$ are comparable to the errors
related to the spatial discretization.  In other words, in the
``optimize-then-discretize'' setting adopted here, gradient
expressions derived based on the continuous formulation,
cf.~\eqref{eq-num-algoS1}, may no longer accurately represent the
sensitivity of the discretized objective function, when the difference
between $u_n$ and $u_{ex}$ is of the order of the space discretization
errors. {This is also confirmed {by computing the errors
    $\|u_n - u_{ex}\|_{H^1}$ at which convergence stagnates.  These
    errors drop} by a factor of roughly 4 when the mesh is refined
  such that $h_{min}$ is reduced by approximately one half (cf.~Figure
  \ref{fig:manuf-exactE}a), as expected from the well-known error
  estimates for the finite-element approximation.}

In Figures \ref{fig:manuf-exactE}a and \ref{fig:manuf-exactE}b it is
evident that the (RCG) method converges much more rapidly (39
iterations) than the (RG) approach (180 iterations).  As expected, the
convergence of the (PG) method (202 iterations) is similar to that of
the (RG) method.  To quantify the convergence rates we use the
following ansatz to represent the errors:
\begin{equation}
\|u_{n}-u_{ex}\|_{2} \sim B_u\, A_u^n, \qquad 
|E_{n}-E_{ex}| \sim B_e\, A_e^n.
\label{eq:manuf-exact-fit}
\end{equation}
The values of the parameters $A_u$ and $A_e$, which represent the
factors by which the corresponding errors are reduced between two
iterations can be obtained from least-squares fits of the data in
Figures \ref{fig:manuf-exactE}a,b in the linear regime. These results
are collected in Table \ref{tab:manuf-exact} (the corresponding fits
are also indicated in Figures \ref{fig:manuf-exactE}a,b). First, these
results demonstrate that the {rate and speed of convergence are}
grid-independent as expected from the general theory of
Sobolev-gradient descent methods \cite{BEC-book-2010-neuberger}. The
data in Table \ref{tab:manuf-exact} can also be interpreted in terms
of the classical theory of the conjugate gradient method {in the
  finite-dimensional Euclidean setting} \cite{nw00}, which for the
minimization of quadratic functions predicts that $A_u \approx
\sqrt{A_e}$. We see that the data from Table \ref{tab:manuf-exact}
satisfies this relationship with the accuracy of a few percent. For
the simple gradient and conjugate gradient methods we furthermore have
the approximate relationships $A_u =(\kappa-1)/(\kappa+1)$ and $A_u
=(\sqrt{\kappa}-1)/(\sqrt{\kappa}+1)$, respectively, where $\kappa$ is
the {``effective''} condition number characterizing the problem.
{It is defined in terms of the condition number of the discrete
  Hessian of the GP energy \eqref{eq-scal-energ} at the minimum
  preconditioned by the metric of the Sobolev space $H_A(\D)$,
  cf.~\eqref{eq-ipHA}, in which optimization is performed.}  Using the
data from Table \ref{tab:manuf-exact} (Mesh 2), we infer that $\kappa
\approx 42.37$ for the gradient (RG) and $\kappa \approx 3.2$ for the
conjugate gradient (RCG) method, indicating that the convergence
acceleration produced in the present problem by the Riemannian
conjugate gradient approach actually exceeds what can be expected from
the standard theory.
\begin{table}[h]
  \begin{center}
    \begin{tabular}{|l|l|l|l|l|l|l|}\hline
      & \multicolumn{3}{|c|}{Mesh 1} & \multicolumn{3}{|c|}{Mesh 2}\\ \hline
      & $A_e$ & $\sqrt{A_e}$ & $A_u$ & $A_e$ & $\sqrt{A_e}$ & $A_u$ \\ \hline
      (RG)     &0.9167&0.9574&0.9496&0.9268&0.9627&0.9538 \\ \hline
      (RCG)	&0.2909&0.5394&0.5275&0.2924&0.5408&0.5238 	\\ \hline	
    \end{tabular}
  \end{center}
  \caption{Parameters characterizing the least-squares fits 
    \eqref{eq:manuf-exact-fit} of the data shown in Figure \ref{fig:manuf-exactE}.}
  \label{tab:manuf-exact}	
\end{table}

\begin{figure}[h]
  \begin{center}
    \includegraphics[width=0.7\textwidth]{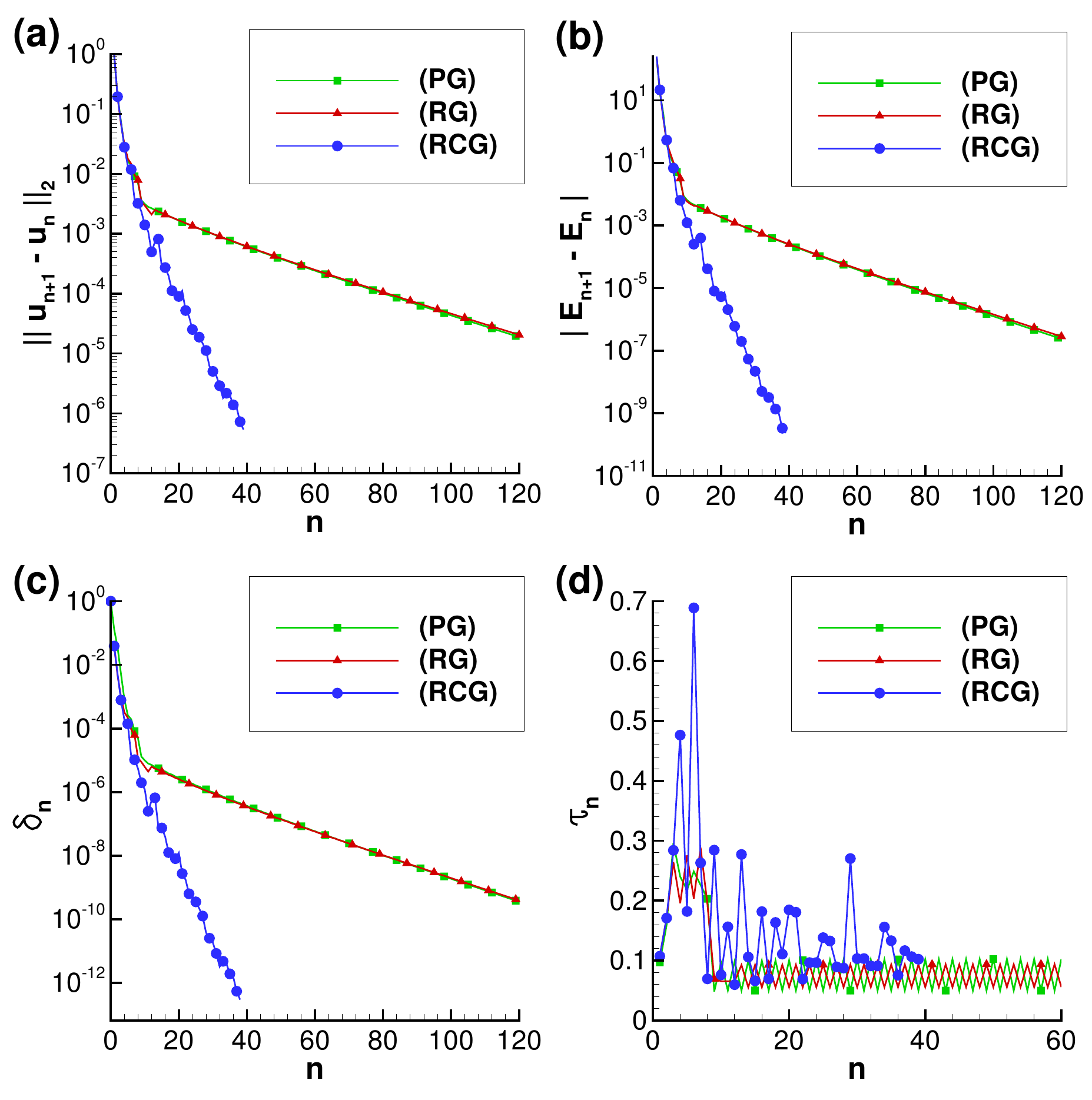}
  \end{center}
  \caption{Test case based on the manufactured solution
    \eqref{eq-manuf-exact}. Convergence of the different quantities
    with iterations $n$: (a) $\|(u_{n+1}-u_n)\|_{2}$, (b)
    $|E_{n+1}-E_n|$, (c) $\delta_n$, cf. \eqref{eq:drift}, drift away
    from the constraint manifold and (d) $\tau_n$, the optimal descent
    step.}
  \label{fig:manuf-exactN}
\end{figure}

Since the exact solution $u_{ex}$ is usually unavailable in physically
relevant problems, we now verify that convergence of iterations can be
monitored based on quantities which do not involve $u_{ex}$. Indeed,
the evolution of $\|u_{n+1}-u_n\|_{2}$ and $|E_{n+1}-E_n|$ with
iterations $n$ shown in Figures \ref{fig:manuf-exactN}a and
\ref{fig:manuf-exactN}b exhibits the same trends as the data shown in
Figures \ref{fig:manuf-exactE}a and \ref{fig:manuf-exactE}b, except
for the slowdown observed in the latter case. This demonstrates that
either of these two quantities can be used to monitor convergence and,
in particular, check the stopping criterion (see also
\cite{BEC-baow-2004-Du,BEC-review-2006-bao,BEC-review-2013-Bao-KRM,BEC-review-2013-antoine-besse-bao,BEC-review-2014-Bao-ICM}).

There are two additional aspects of the convergence of the different
methods we wish to comment on. In Figure \ref{fig:manuf-exactN}c we
show the evolution of the ``drift'' away from the constraint manifold
$\M$ exhibited by the intermediate approximations $\hat{u}_n$ {\em
  before} retraction \eqref{eq:R} is applied
\begin{equation}
  \delta_n =\left| 1 - {\|\hat{u}_n\|^2_{2}}\right|, \qquad n=0,1,\dots.
  \label{eq:drift}
\end{equation}
This quantity measures how far the intermediate steps diverge from the
constraint manifold. We see that, as compared to the (PG) and (RG)
methods, in the (RCG) approach the intermediate approximations always
remain closer to $\M$ .  Finally, the step size $\tau_n$ determined by
the gradient approaches (PG), (RG) and (RCG) via
line/arc-minimization, cf.~\eqref{eq-taun}, \eqref{eq:tau} and
\eqref{eq:tauRCG}, is shown in Figure \ref{fig:manuf-exactN}d. We see
that the step sizes generated by the simple gradient methods, (PG) and
(RG), tend to oscillate between two values, a behavior indicating that
the iterations are trapped in narrow ``valleys''. This is a common
behavior of the steepest descent method when applied to poorly
conditioned problems and is not exhibited by the the (RCG) iterations
where on average the steps also tend to be longer.  {The data in
  Figure \ref{fig:manuf-exactN}a,c,d offers interesting insights about
  the behavior of iterations in different approaches. It follows from
  relation \eqref{eq-norm-error} that $\| u_{n+1} \|_{2}^2 - \|
  u_{n}\|_{2}^2 = \tau_n^2 \| \grdGn \|_{2}^2$ is satisfied for all
  cases considered here, whereas in Figure \ref{fig:manuf-exactN}d it
  is evident that the corresponding step sizes $\tau_n$ are bounded
  away from zero. We therefore deduce that the drift $\delta_n$ is
  directly linked to the magnitude of the gradient $\| \grdGn \|_{2}$,
  which is smallest for the approach with the fastest convergence, \ie
  the RCG method (cf.~Figure \ref{fig:manuf-exactN}a). This
  demonstrates that the small drift away from the constraint manifold
  $\M$ observed in this approach, cf.~Figure \ref{fig:manuf-exactN}c,
  is a consequence of its rapid convergence.}  Performance of
{the} different methods applied to several realistic problems
will be discussed in the next section.

\section{Computation of Rotating Bose-Einstein Condensates}
\label{sec:results}
In this section we compare the performance of the minimization
algorithms (PG), (RG) and (RCG) on a number of test cases involving
configurations of rotating BEC with vortices.  We consider
increasingly complex arrangements: {a single vortex, Abrikosov
  vortex lattices with more than one hundred vortices, {giant
    vortices and, finally, condensates in anisotropic trapping
    potentials}}. To make these test cases more challenging, we
consider large values of the nonlinear interaction constant $C_g$ and
large angular frequencies $C_\Omega$. In some cases, {we will
  provide comparisons between the gradient methods and other
  state-of-the-art techniques, one of which is the (BE) approach
  \eqref{eq-grad-flowBE}--\eqref{eq-steep-normB} implemented using the
  same $P^2$ finite-element setting.  In addition, in order to offer a
  comparison with a higher-order method, we will also solve the
  minimization problem \eqref{eq-Emin} using the library {\tt Ipopt}
  which is interfaced with {\tt FreeFem++}. This approach, which we
  will refer to as (Ipopt), is based on a combination of an interior
  point minimization \cite{Ipopt-int-point-2002}, barrier functions
  \cite{Ipopt-barr-strat-2008} and a filter line-search
  \cite{Ipopt-line-search-2005}. For problems with equality constraints
  only (such as the present problem), (Ipopt) reduces to a Newton-like
  method with an elaborate line-search used to optimally determine the
  step size.} {Here we use (Ipopt) to solve the Euler-Lagrange
  system \eqref{eq-scal-GP-stat} {in the course of which it is
    provided with} the expressions for the $L^2$ gradient and the
  Hessian of the GP energy, reformulated by separating the real and
  imaginary parts of the solution (see also \cite{dan-2016-CPC}).}
{Computations with (Ipopt) are based on several calls of the
  library where, for each call, the residual of the optimality
  condition \eqref{eq-scal-GP-stat}, on which the termination
  criterion is based, is progressively decreased.}

\subsection{Test Case \#1: BEC with a Single Central Vortex}
\label{sec:1v}
\begin{figure}[!h]
  \begin{center}
    \includegraphics[width=0.65\textwidth]{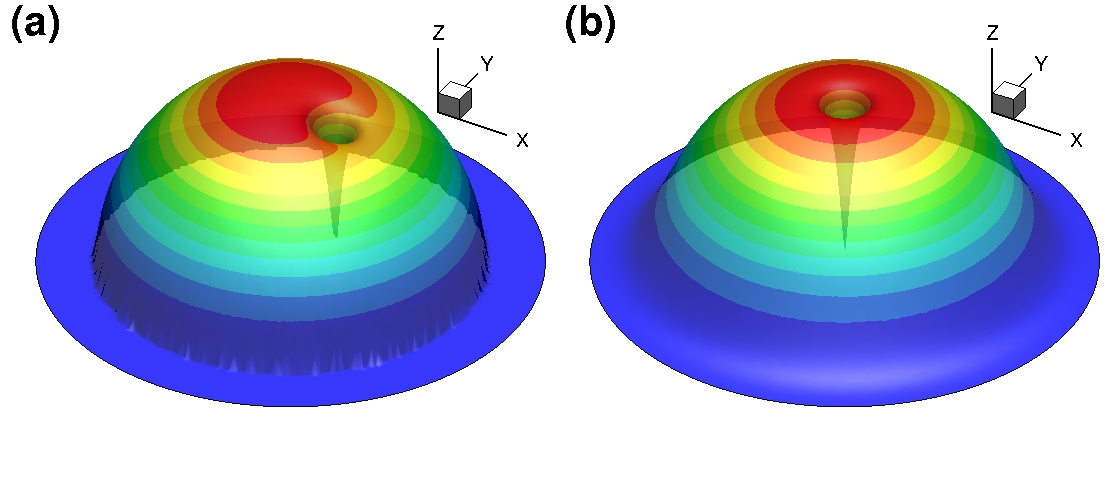}
  \end{center}
  \caption{Computation of a rotating BEC with a single central vortex
    (cf.~\S\ref{sec:1v}). 3D rendering of the atomic density
    $\rho=|u|^2$ for: (a) the initial guess $u_0$ and (b) the
    converged ground state. }
  \label{fig:V1-3D}
\end{figure}

We consider the case of a BEC trapped in a harmonic potential and
rotating at low angular velocities: $C_{\text{trap}}=r^2/2$, $C_g =
500$, $C_\Omega = 0.4$.  For this case, the Thomas-Fermi (TF) theory
\cite{BEC-phys-1999-Stringari-TF} offers a good approximation of the
atomic density $\rho=|u|^2$ of the condensate $\rho \approx \rtf =
\left((\mu - C^\text{eff}_\text{trap}) / {C_g} \right)_+$ with the
effective trapping potential $C^\text{eff}_\text{trap}$ given by
\eqref{eq-scal-Ctrap-eff}.  By imposing $\int_{\cal D} \rtf = 1$, we
can derive analytical expressions for the corresponding approximation
of the chemical potential {$\mu \in \RR$} \cite{dan-2016-CPC}.
From the same TF approximation, we can estimate the radius of the
condensate as $R_\TF= \sqrt{2\mu/(1-C_\Omega)}$.  Consequently, we set
up the computational domain $\D$ as a disk of radius $R= 1.25 R_\TF=
6.56$. The initial guess $u_0$ is taken in the form of an off-center
vortex placed at $(x_v=0.25, y_v=0)$, see Figure \ref{fig:V1-3D}a.  We
use the ansatz $u_0 = \sqrt{\rtf} \, u_v$, where $u_v = r /
\sqrt{r^2+2\xi^2}\, e^{i \theta}$ with $(r,\theta)$ representing the
polar coordinates centered at $(x_v,y_v)$ and $\xi = 1/\sqrt{2 \mu}$
the non-dimensional healing length, which is a good approximation of
the vortex radius in rotating BEC \cite{BEC-physV-2001-fetter}.
{A similar} ansatz will also be used in subsequent sections to
set up initial guesses with vortices for the calculation of more
complicated BEC configurations. The stopping criterion
\eqref{eq:stopping} is used with the value
$\varepsilon_{st}=10^{-12}$. In these calculations the grid remains
fixed (\ie no grid adaptation is performed), {with 9,578 vertices
  and 18,825 triangles}. In the (BE) method the imaginary time step is
chosen as $\delta t = 10$, which proved to be optimal for convergence
after testing values of $\delta t$ in the range from $0.01$ to $100$.
{In the (Ipopt) approach, two successive calls to the library
  were sufficient {to converge the solution to the same level of
    accuracy as with other methods}.}

For the considered physical parameters the ground state features a
vortex centered at the origin.  All considered methods, (PG), (RG),
(RCG), {(BE) and (Ipopt)}, converged to the same ground state
shown in Figure \ref{fig:V1-3D}b. In order to assess their respective
rates of convergence, we compute a reference (``exact'') solution
$u_{ex}$ using the same grid and starting the minimization algorithms
from the initial guess $u_0$ with $(x_v=0, y_v=0)$. The corresponding
energy and angular momentum are $E_{ex}=8.36059$ and $L_{ex}=1$.

\begin{figure}[!h]
  \begin{center}
    \includegraphics[width=0.6\textwidth]{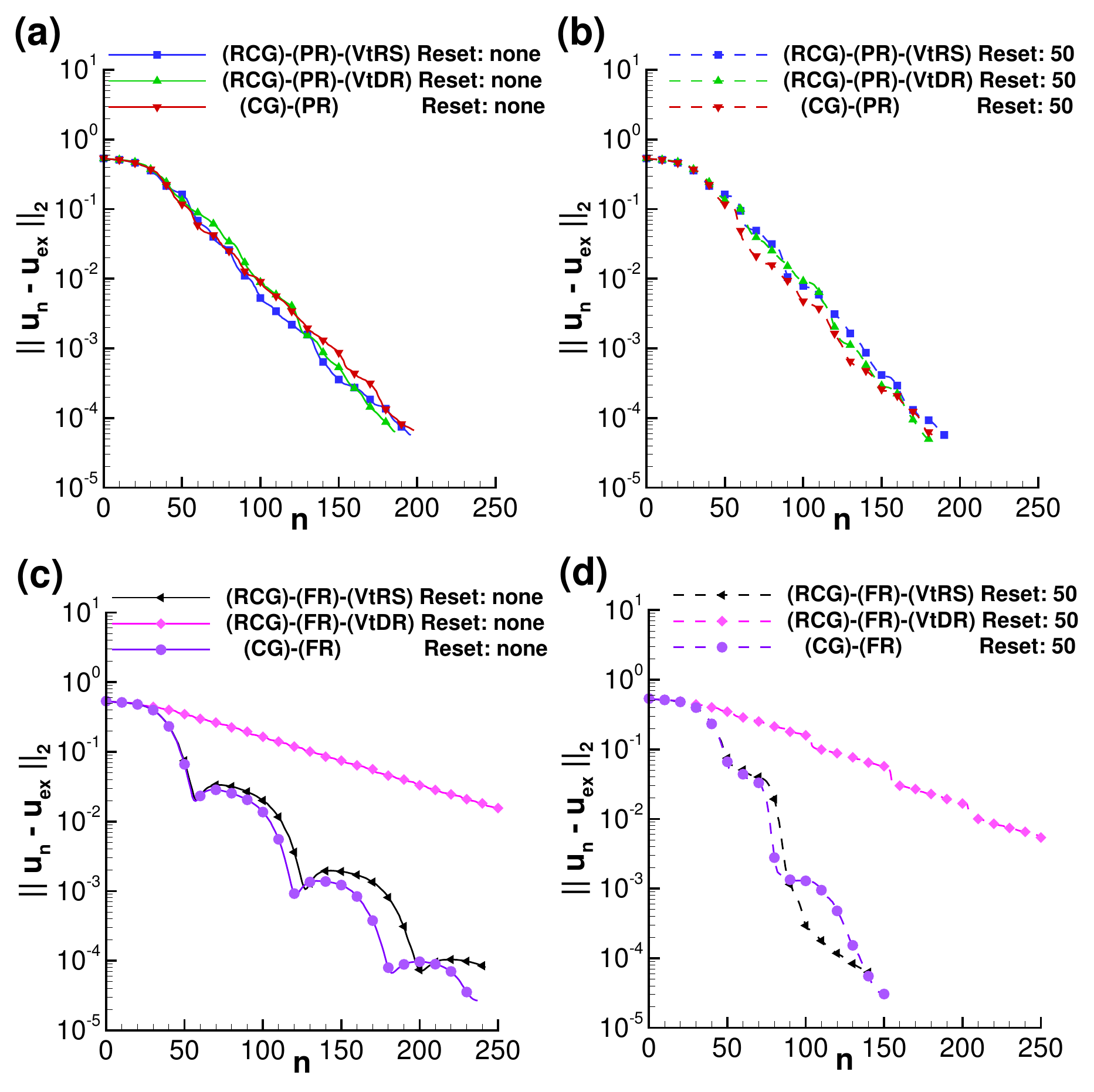}
  \end{center}
  \caption{{Computation of a rotating BEC with a single central vortex
    (cf.~\S\ref{sec:1v}). Comparison of different design choices for
    the (RCG) method in terms of convergence of {the error norm}
    $\|u_{n}-u_{ex}\|_{2}$.}}
  \label{fig:V1-conv-num-RCG}
\end{figure}

{The performance of the approaches corresponding to the
  different design choices discussed in \S \ref{sec:design} is
  summarized in Table \ref{tab:V1-perf}, whereas their convergence
  monitored in terms the error norm $\|u_{n}-u_{ex}\|_{2}$ is
  illustrated in Figure \ref{fig:V1-conv-num-RCG}. As compared with
  the results analyzed in \S \ref{sec:res-manuf}, the difference here
  is that we now allowed for periodic resets of the momentum term
  $\beta_n$ to zero and it was found by trial-and-error that the
  fastest convergence in terms of the CPU time was obtained when such
  resets were performed every 50 iterations. From Table
  \ref{tab:V1-perf} we conclude that the (RCG)-(PR) approach is more
  robust than the (RCG)-(FR) approach with respect to the choice of
  the vector transport and the reset frequency.  We also note that
  ignoring the vector transport, which is the case in the (CG)
  methods, produces a significant increase of the computational time
  resulting from slow convergence of the arc-minimization procedure.
  As regards resetting the momentum term to zero, we note that it
  turns out to be particularly important for the (FR) approaches,
  concurring with the insights about them already known from the
  Euclidean setting \cite{nw00}. In addition to being costly to
  determine, optimal momentum reset strategies also tend to be
  strongly problem-dependent. Thus, with this in mind, we conclude
  that the (RCG)-(PR)-(VtRS) approach again turns out to be the most
  efficient and is also characterized by the most consistent
  performance. Therefore, it will be used for further tests in the
  remainder of this section and we will continue to refer to is as
  (RCG).  }

\begin{table}[h]
  \begin{center}
    \begin{tabular}{|r||r|r|r|r|}\hline
      Method & \multicolumn{2}{|c|}{Reset: none} & \multicolumn{2}{|c|}{Reset: 50  iter}\\ \hline
             & iter & CPU & iter & CPU\\ \hline \hline
      (RCG)-(PR)-(VtRS)  &196 &913 &190 &903\\ \hline
      (RCG)-(PR)-(VtDR) & 186 &954 &180 &868\\ \hline
      (CG)-(PR)             & 198 &1328 &184 &1271\\ \hline \hline	
      (RCG)-(FR)-(VtRS)  & 242 &1244 &140 &637\\ \hline
      (RCG)-(FR)-(VtDR) &  563 &2652 &359 &1703\\ \hline
      (CG)-(FR)             & 237 &1583 &150 &1077\\ \hline \hline	
      (RG)  &2643 & 13872 & &\\ \hline
      (PG)  & 3631 & 8472 &  &\\ \hline
      (BE)  &2796 & 6838 &    &\\ \hline		
      (Ipopt)  &18 & 99  &    & \\ \hline			
    \end{tabular}
%
  \end{center}
	
  \caption{{Computation of a rotating BEC with a single central vortex
    (cf.~\S\ref{sec:1v}).  Performance of {the gradient methods
      corresponding to the different design choices, cf.~\S
      \ref{sec:design}, and of the (Ipopt) approach measured} in terms
    of the number of iterations (iter) and the computational time in
    seconds (CPU) {required for} convergence. }}
  \label{tab:V1-perf}	
\end{table}

\begin{figure}[!h]
  \begin{center}
    \includegraphics[width=0.6\textwidth]{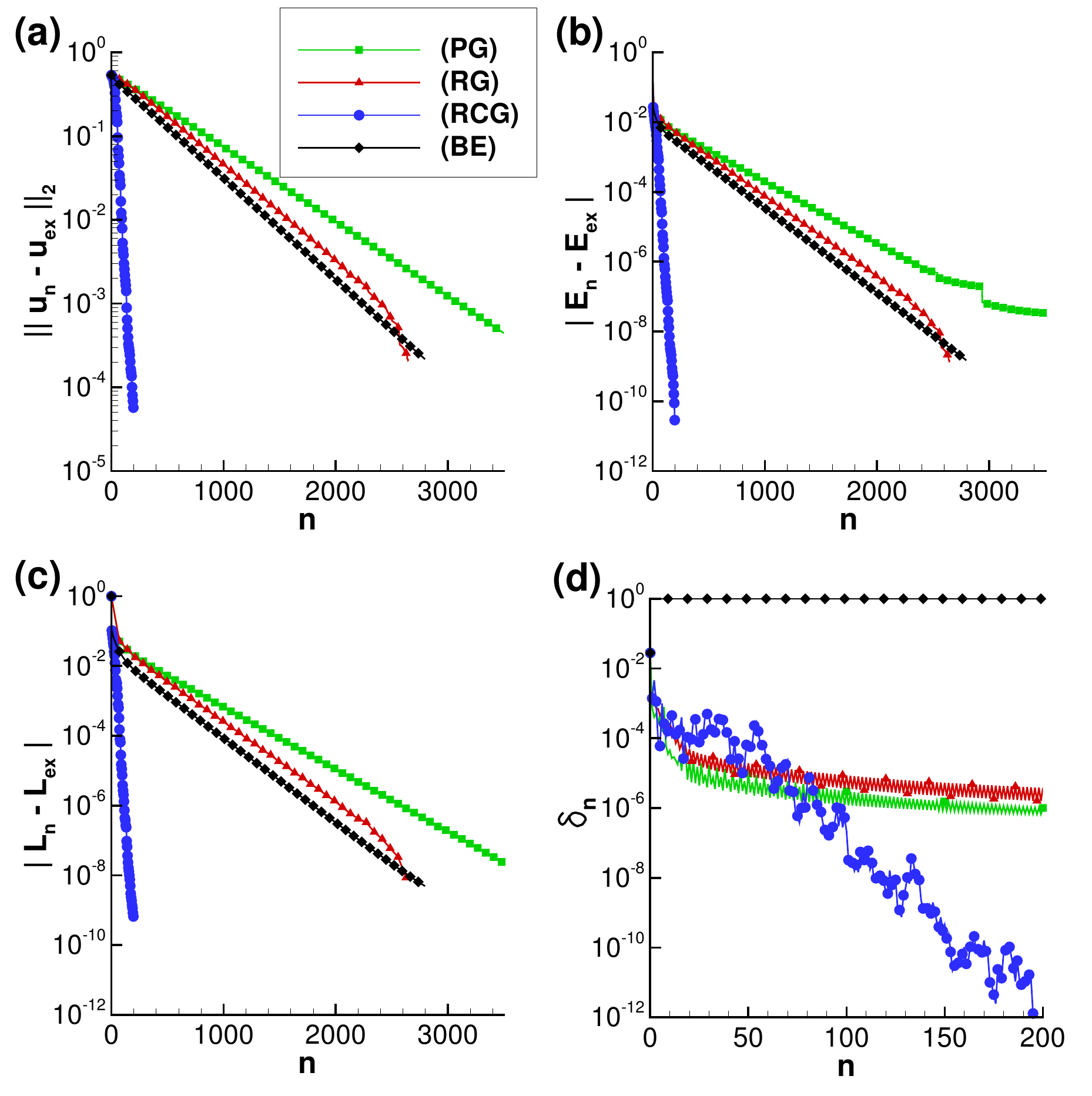}
  \end{center}
  \caption{Computation of a rotating BEC with a single central vortex
    (cf.~\S\ref{sec:1v}). Convergence of the different quantities with
    iterations $n$: (a) $\|(u_{n}-u_{ex})\|_{2}$,
   (b) $|E_{n}-E_{ex}|$, (c)  
    $|L_{n}-L_{ex}|$ and (d) $\delta_n$, cf.~\eqref{eq:drift}, the drift away from the manifold.}
  \label{fig:V1-conv-num}
\end{figure}

The quantities $\|u_n-u_{ex}\|_{2}$, $|E_n-E_{ex}|$ and
$|L_{n}-L_{ex}|$ shown in Figures \ref{fig:V1-conv-num}a,
\ref{fig:V1-conv-num}b and \ref{fig:V1-conv-num}c as functions of $n$
indicate that while the (RG) and (BE) methods converge with similar
rates, the (RCG) approach converges much faster.  The drift $\delta_n$
away from the constraint manifold $\M$ at intermediate steps
$\hat{u}_n$, cf.~\eqref{eq:drift}, during the first 200 iterations is
shown for different methods in Figure \ref{fig:V1-conv-num}d. We note
that {for the (BE) method this quantity is always
  $\mathcal{O}(1)$ which is due to the fact that the RHS of
  \eqref{eq-grad-flowBE} is based on an unprojected gradient further
  compounded by a large step size $\delta t$ used.}  The normalization
step is therefore crucial in this approach.  On the other hand,
$\delta_n$ is reduced faster in the (RCG) approach where it also
attains lower values than in the (PG) and (RG) methods.

{Finally, let us note that, even though the (RCG) method
  outperforms all other first-order methods, the (Ipopt) approach
  actually converges much faster, cf.~Table \ref{tab:V1-perf}. Only 18
  Hessian evaluations are needed during the two calls to {\tt Ipopt}
  and the total computational time is {smaller} by a factor of
  10 as compared to the (RCG) methods. {However, we stress} that
  this test was performed with a fixed grid {and in fact} this
  remarkable performance of the (Ipopt) method will be lost when
  computing more complicated cases that require mesh adaptivity (see
  the next section).}

\subsection{Test Case \#2: BEC with a dense Abrikosov vortex lattice}
\label{sec:AL}

We now move on to consider more challenging test cases corresponding
to a harmonic potential ($C_\text{trap}=r^2/2$), high rotation rate
($C_\Omega=0.9$) and large values of the nonlinear interaction
constant ($C_g=1000$ to $15 000$). We note that for the harmonic
trapping potential there is a physical limit occurring at the rotation
frequency $C_\Omega=1$ when the trapping is canceled by the
centrifugal force (\ie $C^\text{eff}_\text{trap}=0$, see
eq.~\eqref{eq-scal-Ctrap-eff}). The next section will consider cases
with a modified trapping potential, allowing for higher rotation
frequencies.

\begin{figure}[h]
  \begin{center}
    \includegraphics[width=0.6\textwidth]{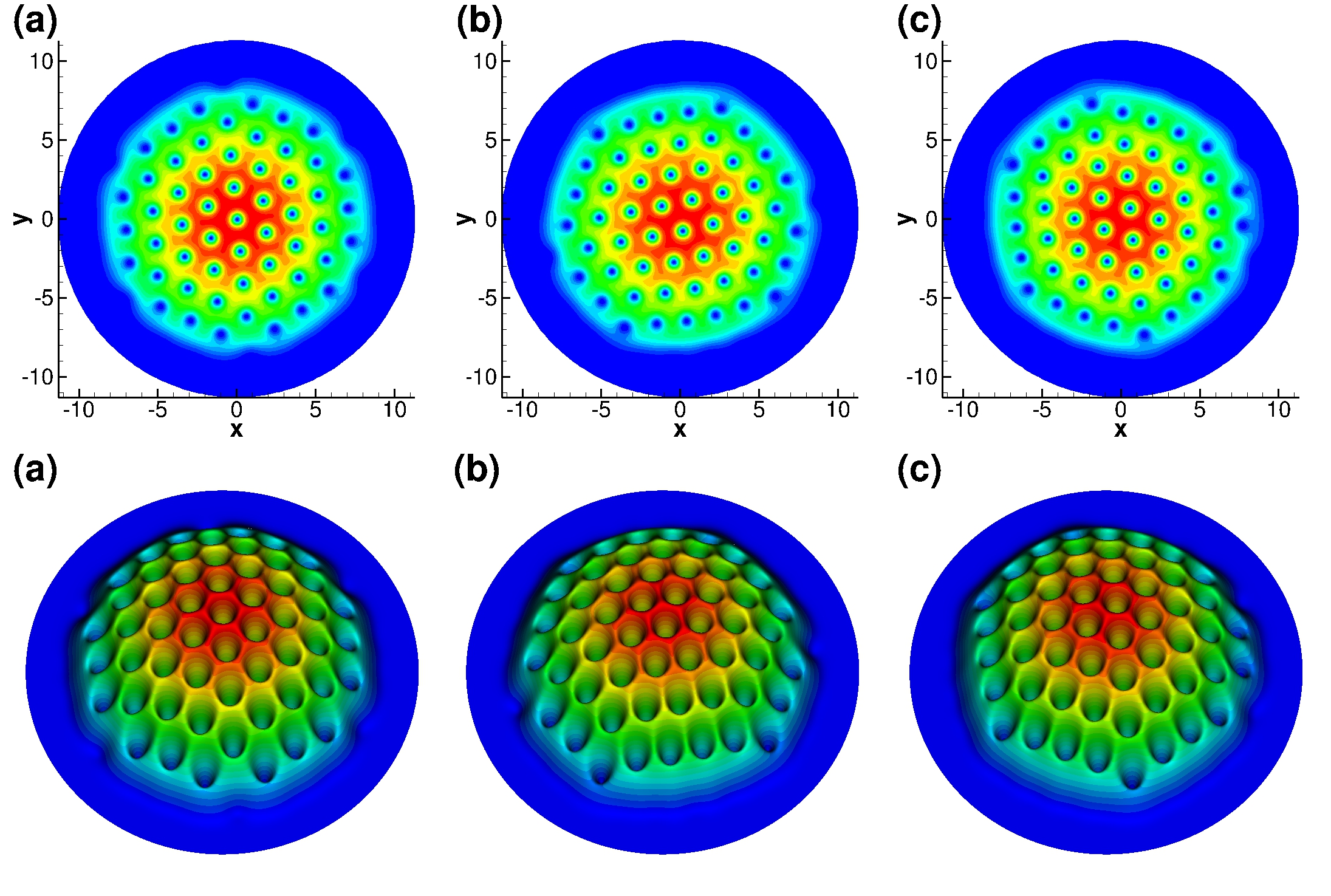}
  \end{center}
  \caption{Computation of a rotating BEC with a dense Abrikosov vortex
    lattice (cf.~\S\ref{sec:AL}). Stationary states obtained using
    {the (RCG)-(PR)-(VtRS) method with mesh adaptivity} and
    different initial conditions $u_0$: (a) ``ansatz d'' suggested in
    \cite{BEC-baow-2016-Newton}, (b) Thomas-Fermi atomic density with
    one central vortex, (c) Thomas-Fermi atomic density with a ring of
    six vortices. The figures in the first row show contours of the
    atomic density (normalized by its maximum value $\rho/\rho_{max}$)
    and in the second row they show the 3D-rendering of the same
    contours.  The corresponding energies are: (a) $E=6.3615$, (b)
    $E=6.3621$, (c) $E=6.3620$, to be compared to the reference value
    $E^{ref} = 6.3607$ from \cite{BEC-baow-2016-Newton}.}
  \label{fig:AL-isos}
\end{figure}

We start with the case with $C_g=1000$ and $C_\Omega=0.9$ for which
the ground state features over fifty vortices arranged in a regular
triangular lattice called the Abrikosov lattice.  The difficulty here
is to obtain a very regular lattice, in particular for the vortices
located near the border of the condensate where the atomic density is
low. This explains why the results previously reported for this case
exhibit configurations with somewhat different arrangements of the
peripheral vortices which nevertheless have very similar energy levels
\cite{BEC-CPCm-2009-zeng,dan-2010-SISC,BEC-numm-2014-antoine-duboscq-JCP}.
These differences can be attributed to the use of different initial
guesses $u_0$.  A nearly perfect arrangement of vortices on a
triangular/hexagonal lattice is reported in the recent study
\cite{BEC-baow-2016-Newton} and will be considered here as a reference
result used to validate our methods (the corresponding energy level is
$E^{ref} = 6.3607$). Details of the computed stationary states depend
on the initial guess $u_0$, and we used three distinct forms of $u_0$:
(i) ``ansatz d'' proposed in \cite{BEC-baow-2016-Newton} to model a
central vortex using Gaussian functions, and the Thomas-Fermi
approximation described in \S\ref{sec:1v} with (ii) one central vortex
and (iii) six vortices. The corresponding stationary solutions
{obtained using the (RCG)-(PR)-(VtRS) method are shown in Figures
  \ref{fig:AL-isos}a, \ref{fig:AL-isos}b and \ref{fig:AL-isos}c.  We}
can see that the central parts of the vortex lattices are in all cases
essentially identical (modulo rotation) and some differences are
detected among the peripheral vortices. The values of energy
corresponding to these configurations differ by less than $0.01\%$.

In the following we carry out computations starting from the initial
guess (i) and {mesh} adaptation is now performed during
iterations {which are declared converged when the termination
  condition \eqref{eq:stopping} with $\varepsilon_{st}=10^{-9}$ is
  met.  For this challenging test case the performance of a few
  selected design choices for the (RCG) and (CG) approaches is
  summarized in Table \ref{tab:AL-perf}. Since the (FR) methods with
  and without momentum resets failed to converge to the same minimum
  as other approaches, we focus here on the (PR) techniques and note
  the fast convergence of the (RCG)-(PR)-(VtRS) approach which will be
  used in further tests in this section; for brevity, we will refer to
  it simply as (RCG).  It is interesting to note in Table
  \ref{tab:AL-perf} that the performance of the (Ipopt) method is
  significantly degraded as compared to the results from \S
  \ref{sec:1v} and is now comparable to that of the (RCG) method.  The
  reason is that {\tt Ipopt} is linked as an external library to {\tt
    FreeFem++} and therefore we cannot directly use mesh adaptivity in
  its internal algorithm. As a result, one has to use an external
  algorithm to couple the computation of the minimizer with the mesh
  adaptivity procedure employed in the other methods.  The
  computations in the present case required 17 calls to {\tt Ipopt},
  with a total of 354 Hessian evaluations and a large number of
  internal iterations.}

\begin{table}[h]
  \begin{center}
    \begin{tabular}{|r||r|r|l|}\hline
      {Method}  & $E$& iter & CPU\\ \hline \hline
      (RCG)-(PR)-(VtRS)&6.3615 &1339& 23192\\ \hline
      (RCG)-(PR)-(VtDR)&6.3615 &2684& 59431\\ \hline
      (CG)-(PR)            &6.3615 &1557& 25693\\ \hline
      (Ipopt)&6.3621 &354& 22943\\ \hline						
    \end{tabular}
    %
    
  \end{center}
  
  \caption{{Computation of a rotating BEC with a dense Abrikosov vortex
    using mesh adaptivity and initial condition $u_0$ described by the
    ``ansatz d'' suggested in \cite{BEC-baow-2016-Newton}
    (cf.~\S\ref{sec:AL}). Performance of {the gradient methods
      corresponding to the indicated design choices, cf.~\S
      \ref{sec:design}, and of the (Ipopt) approach measured} in terms
    of the number of iterations (iter) and the computational time in
    seconds (CPU) {required for} convergence. }}
  \label{tab:AL-perf}	
\end{table}

\begin{figure}[h]
	\begin{center}
		\includegraphics[width=\textwidth]{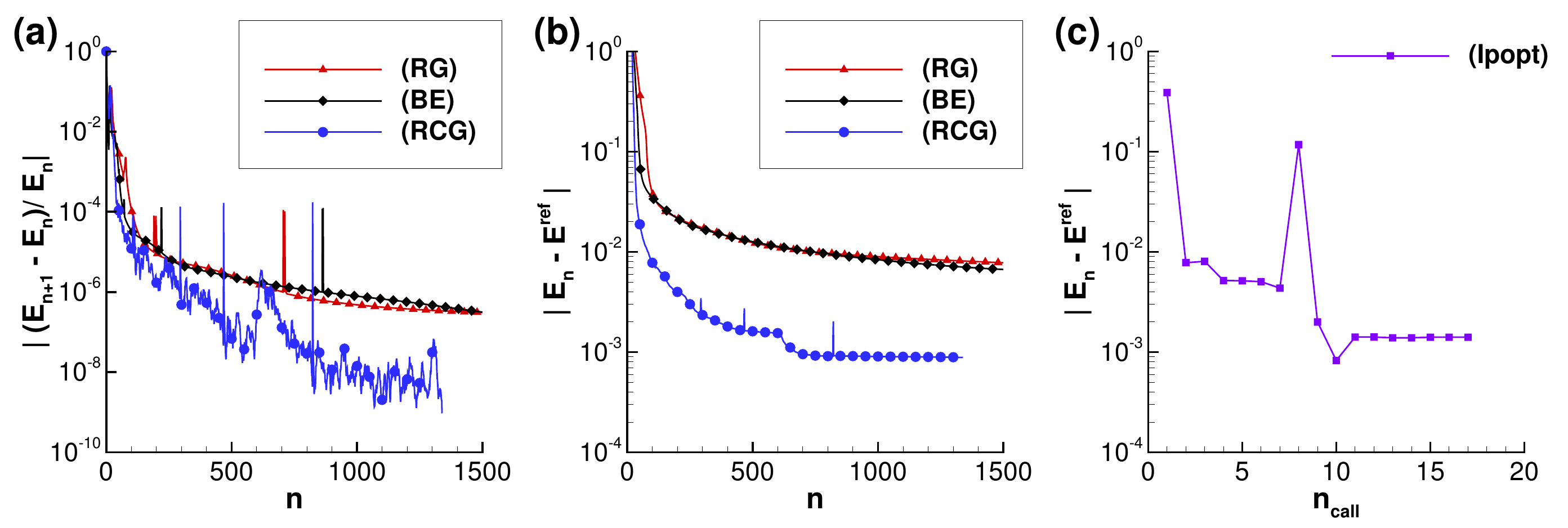}
	\end{center}
	\caption{Computation of a rotating BEC with a dense Abrikosov vortex
		lattice (cf.~\S\ref{sec:AL}). Convergence of the different
		quantities with iterations $n$: (a) $|(E_{n+1}-E_n)/E_n|$ and
		{(b, c)} $|E_{n}-E^{ref}|$. All computations start from the
		initial guess suggested in \cite{BEC-baow-2016-Newton} ("ansatz
		d").}
	\label{fig:AL-conv-num}
\end{figure}

{Convergence of the iterations carried out with the (RG), (RCG)
  and (BE) methods is compared in Figures \ref{fig:AL-conv-num}a and
  Figures \ref{fig:AL-conv-num}b.}

 {The gradient (RG) and the
  backward-Euler (BE) methods show a similar, but markedly slower
  convergence (they were stopped after 5000 iterations) than the (RCG)
  method.}  We note that the peaks in the curves shown in Figures
\ref{fig:AL-conv-num}a and \ref{fig:AL-conv-num}b result from
reinterpolation of intermediate solutions after grid adaptation
\cite{dan-2010-JCP,dan-2016-CPC}. In these figures we see that in all
cases convergence slows down at later iterations which is related to
the slow rearrangement of vortices near the boundary of the
condensate. Another possible reason is that since in
\cite{BEC-baow-2016-Newton} a different discretization was used
(Fourier spectral approach with periodic boundary conditions), the
value of $E^{ref}$ taken from that reference might not exactly
correspond to ours.  {Since convergence is monitored differently
  for the (Ipopt) method, in Figure \ref{fig:AL-conv-num}c we show the
  decrease of the energy with the number of calls to the library
  {(we add that the energy also tends to exhibit significant
    oscillations during iterations performed within each such call).
    While the energy level attained with the (Ipopt) method is a
    little higher than obtained with the (RCG) and (CG) methods,
    the corresponding minimizer is similar to that shown in Figure
    \ref{fig:AL-isos}c.}}

\subsection{Test Case \#3: BEC with a large number of vortices}
\label{sec:AL2}

\begin{figure}[h]
  \begin{center}
    \includegraphics[width=0.5\textwidth]{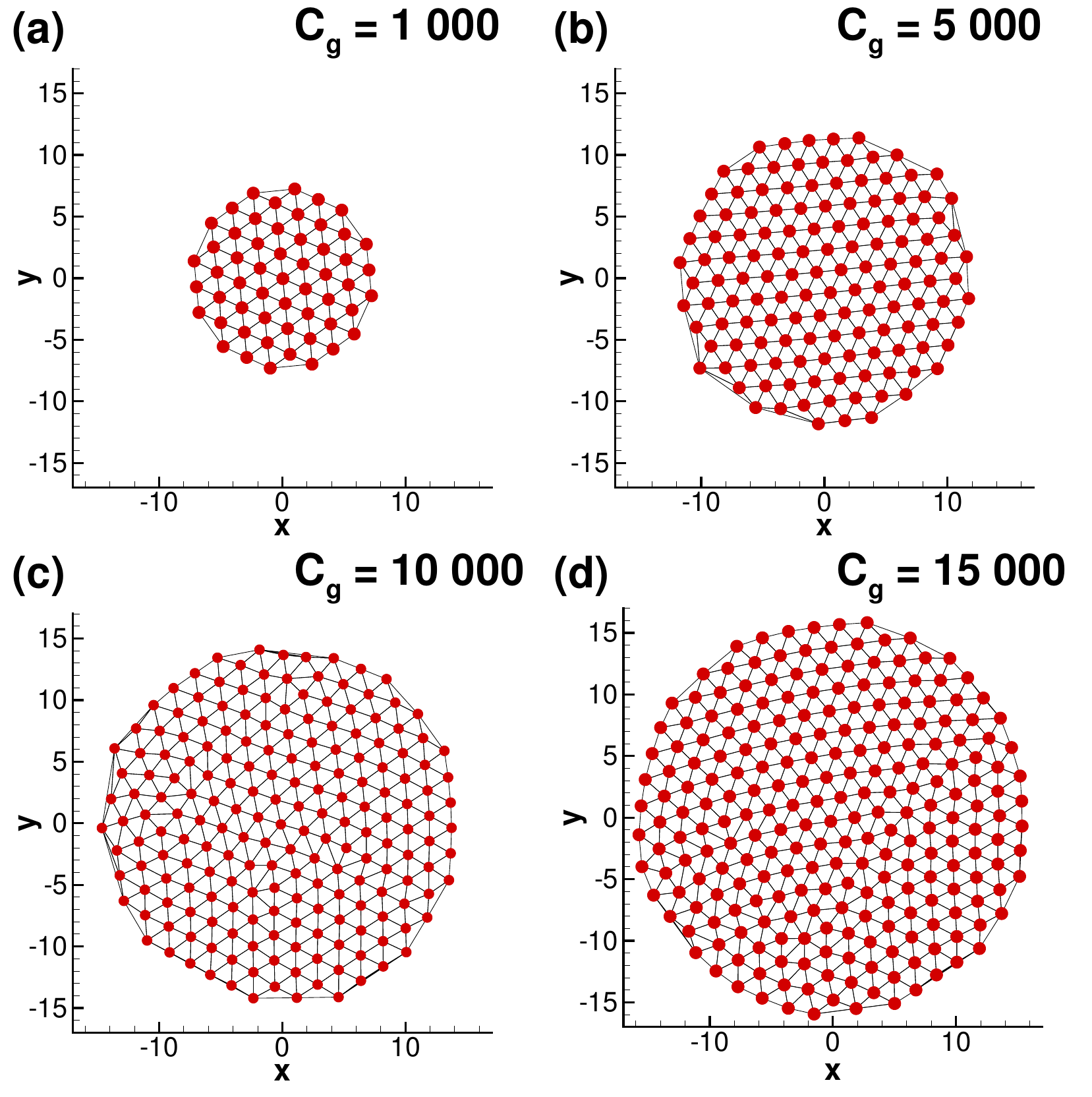}
  \end{center}
  \caption{Computation of a fast rotating BEC ($C_\Omega=0.9$) in a
    harmonic trapping potential (cf.~\S\ref{sec:AL}).  The Abrikosov
    vortex lattice is represented using the Delaunay triangulation
    built from the detected vortex centers. Configurations obtained
    for large values of the nonlinear interaction constant: $C_g=1
    000$ (55 vortices), $C_g=5 000$ (134 vortices), $C_g=10 000$ (193
    vortices) and $C_g=15 000$ (237 vortices).}
  \label{fig:AL-5k}
\end{figure}

{In this section, we use the (RCG)-(PR)-(VtRS) method} to compute
fast rotating BEC ($C_\Omega=0.9$) corresponding to large values of
the nonlinear interaction constant with $C_g$ varying from 1,000 to
15,000.  For these difficult cases, a more physically relevant
assessment of the convergence of iterations is provided by the
alignment of vortices on parallel lines inside the vortex lattice.
Since isocontours of atomic density do not always coincide with these
lines, we developed a post-processing approach to identify the centers
of vortices by detecting local minima of the function $\rtf-\rho$.
This post-processing is similar to that used for experimental data
\cite{BEC-physV-2001-codd} or 3D numerical simulations \cite{dan-2005}
and also allows to build the Delaunay triangulation of the lattice and
compute the radius of each vortex.  The resulting stationary states
are presented in this way in Figure \ref{fig:AL-5k}. We notice an
arrangement of vortices on a nearly perfect lattice for $C_g=1,000$
and $5,000$, and a less regular arrangement for $C_g=10,000$ and
$15,000$ with the presence of some defects in the lattice. This effect
could be related to physical theories addressing the non-uniformity of
the inter-vortex spacing in dense Abrikosov lattices
\cite{BEC-physV-2001-codd,BEC-physV-2004-sheehy}.

\subsection{Test Case \#4: BEC with giant vortex}
\label{sec:GIANT}

To overcome the limit $C_\Omega=1$ imposed by the harmonic trapping
potential, a modified ``harmonic-plus-Gaussian'' potential was tested
in experiments \cite{BEC-physV-2004-bretin}. In
\cite{dan-2004-aft,dan-2005} this new experimental set-up was modeled
as
\begin{equation}
\label{eq-scal-trap-Q}
C_\text{trap}(x,{y})=(1-\alpha) r^2 + \frac{1}{4} k r^4, 
\end{equation}
with the possibility to switch from a ``quartic-plus-quadratic''
potential ($\alpha < 1$), which corresponds to experiments, to a
``quartic-minus-quadratic'' potential ($\alpha > 1$), which is
experimentally feasible but was never tested. Adapting the analysis
from \cite{dan-2004-aft} to our 2D case, we obtain three possible
regimes depending on the type of potential:\\
$\bullet$ Regime 1: ``quartic-plus-quadratic'' (or weak attractive)
potential obtained for $\alpha < 1$ and $\mu > 0$ (see the
Thomas-Fermi approximation in \S\ref{sec:1v}).\\
$\bullet$ Regime 2: weak ``quartic-minus-quadratic'' (or weak
repulsive) potential obtained when $\alpha > 1$ and $\mu > 0$; this
regime appears when
$|1-\alpha| < \left( k \sqrt{3 C_g/\pi}\right)/2$.\\
$\bullet$ Regime 3: strong ``quartic-minus-quadratic'' (or strong
repulsive) potential obtained when $\alpha > 1$ and $\mu < 0$; this
regime appears when $|1-\alpha| > \left( k \sqrt{3 C_g/\pi}\right)/2$.

\begin{figure}
  \begin{center}
    \includegraphics[width=0.55\textwidth]{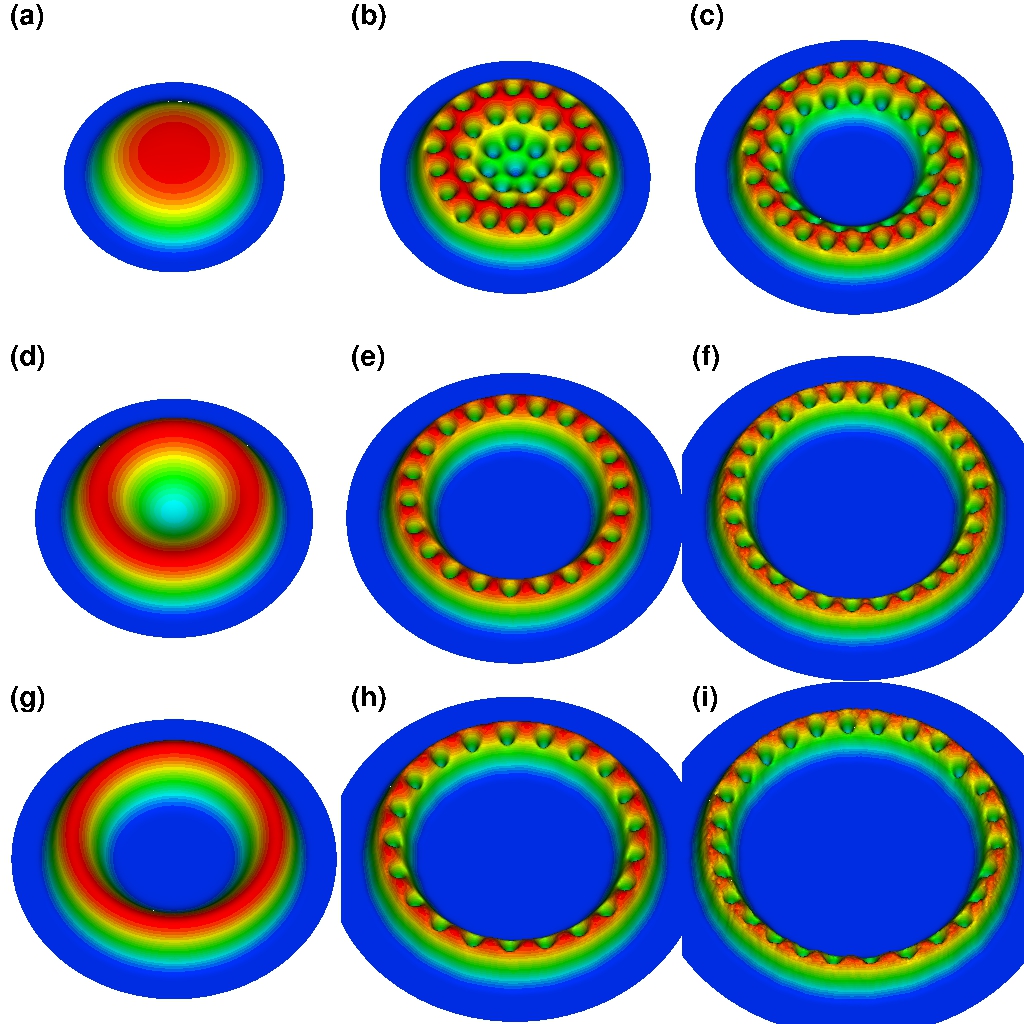}
  \end{center}
  \caption{Computation of a rotating BEC with giant vortices
    (cf.~\S\ref{sec:GIANT}). 3D-rendering of the atomic density (normalized by its maximum
    value $\rho/\rho_{max}$) obtained in Regime 1
    (a, b, c), Regime 2 (d, e, f) and Regime 3 (g, h, i) for
    different rotations $C_\Omega=0$ (first column), $C_\Omega=3$
    (second column) and $C_\Omega=4$ (third column).}
  \label{fig:VG-isos}
\end{figure}
All computations are performed with the {(RGC)-(PR)-(VrRS)}
method and the obtained stationary states are presented in Figure
\ref{fig:VG-isos}. The parameters for these simulations are:
$C_g=1000$, $k=1$ and $\alpha=1/2$ (Regime 1), $\alpha=11/2$ (Regime
2), $\alpha=9$ (Regime 3).  In the first column of Figure
\ref{fig:VG-isos} we notice that the atomic density distribution in
the condensate without rotation ($C_\Omega=0$) changes from the
classical parabolic profile in Regime 1 to a Mexican-hat type profile
in Regime 2 and, finally, to a profile with a central hole in Regime
3. It is then expected that, when rotation is applied, in Regimes 2
and 3 the condensate will develop a central hole (or giant vortex) at
lower rotations frequencies than in Regime 1. This prediction is
indeed supported by the results in Figure \ref{fig:VG-isos} (second
and third column). When rotation is increased, the condensate
configuration evolves from a classical (Abrikosov) vortex lattice to a
vortex lattice with a central depletion and, finally, to a giant
vortex surrounded by a ring of individual vortices. The giant vortex
is indeed obtained for lowest rotation frequencies in Regimes 2 and 3.
The existence of a giant vortex was predicted theoretically (\eg
\cite{BEC-physV-2001-fetter}) and was already observed in 2D
\cite{BEC-physV-2002-kasamatsu-giant} and 3D computations
\cite{dan-2004-aft,dan-2005}.

\subsection{Test Case \#5:  Strongly Anisotropic BEC with Vortices}
\label{sec:ANISO}
\begin{figure}
  \begin{center}
    \includegraphics[width=0.75\textwidth]{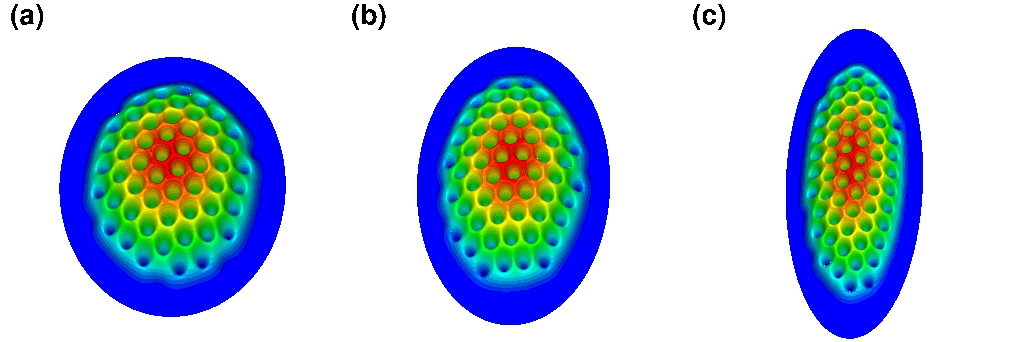}
  \end{center}
  \caption{Computation of a rotating BEC in a strongly anisotropic
    trapping potential (cf.~\S\ref{sec:ANISO}).   3D-rendering of the atomic density (normalized by its maximum
    value $\rho/\rho_{max}$) obtained for
    $\epsilon=0.15$ (a),
    $\epsilon=0.35$ (b), $\epsilon=0.65$ (c).  }
  \label{fig:ANISO-isos}
\end{figure}

To {demonstrate} the efficiency of {the (RGC)-(PR)-(VrRS)}
method for the case without any symmetry, in this section we consider
a rotating BEC trapped in a strongly anisotropic (asymmetrical)
potential of the form suggested in \cite{BEC-physV-2004-ANISO}
\begin{equation} 
\label{eq-scal-trap-ANISO}
C_\text{trap}(x,{y})=\frac{1}{2}\left[(1+\eta^2) x^2  + (1-\eta) y^2\right], \quad \eta=2(1-C_\Omega)\epsilon,
\end{equation}
where $\epsilon < 1$ characterizes the anisotropy of the trap for very
high rotation frequencies $C_\Omega \approx 1$.  The theoretical
analysis presented in \cite{BEC-physV-2004-ANISO} shows that, when the
condensate contains a large number of vortices, the deviation of the
vortex lattice from a triangular arrangement is small. This finding is
supported by our computational results shown in Figure
\ref{fig:ANISO-isos} for three values of the anisotropy parameter
$\epsilon$. This example illustrates the flexibility of the
finite-element discretization, cf.~\S\ref{sec:numer}, in handling
highly deformed computational domains $\D$.

\section{Conclusions}
\label{sec:final}
%
The difficulty of direct minimization of the Gross-Pitaevskii energy
functional with rotation comes from the unit-norm constraint
\eqref{eq-scal-cons}.  The novel idea proposed here is to transform
this problem to an unconstrained Riemannian optimization problem
defined on a spherical manifold and then develop a Riemannian
conjugate gradient (RCG) method based on classical approaches.  The
key {ingredients} of this new method are the following: (i) the
gradient direction is derived using the theory of Sobolev gradients
and relies on a physically-inspired definition of the inner product
which accounts for rotation \cite{dan-2010-SISC}, thereby offering a
good preconditioning for the problem; (ii) the gradient is projected
on the subspace tangent to the spherical manifold before being used in
simple gradient or conjugate gradient methods, which ensures the
iterates stay close to manifold $\M$; (iii) the conjugate descent
direction is computed using classical approaches (\ie the
Polak-Ribi{\`{e}}re or Fletcher-Reeves variant of the nonlinear
conjugate gradient method) and the Riemannian vector transport is used
to bring the gradient and descent directions determined at the
previous iteration to the current tangent subspace $\TunM$; (iv) the
optimal descent step is computed by solving an arc-minimization
problem (in which samples are constrained to lie on the manifold),
instead of the classical line-minimization; (v) finally, the updated
solution is ``retracted'' back to the spherical manifold. {In
  our study we carefully analyzed the effect of the key design
  choices, namely, the form of the momentum term and of the vector
  transport, on the performance of the Riemannian conjugate gradients
  approach. Based on tests involving several different problems, we
  conclude that the (RGC)-(PR)-(VrRS) approach, combining the
  Polak-Ribi{\`{e}}re form of the momentum term with the vector
  transport based on the Riemannian submanifold structure, exhibited
  the most robust and efficient performance. The (RCG)-(PR) methods in
  general also showed a systematic improvement over the (CG)
  approaches without vector transport.}

As demonstrated by our tests performed in the finite-element setting,
several features make the (RCG) method very appealing for practical
computations: (i) since the ``optimize-then-discretize'' paradigm is
used, the preconditioning is mesh-independent; (ii) the Riemannian
retraction and transport operators are simple to implement; (iii) for
the arc-minimization problem a classical approach such as Brent's
method can be easily adapted; (iv) there are no tuning parameters or
trust-region tests involved.  In addition, general mesh refinement or
mesh adaptivity strategies are compatible with the RCG method without
any modifications. Our extensive numerical experiments showed a
significant improvement of the convergence rate of the RCG method over
the simple gradient and imaginary-time methods. {For more
  involved problems requiring mesh adaptation the (RCG) approach
  exhibited performance comparable to the (Ipopt) method which for
  equality-constrained problems implements a Newton-type technique.
  The reason is that there is no straightforward way to incorporate
  mesh adaptation in the (Ipopt) approach, something that can be done
  rather easily in the (RCG) method. We stress that the use of mesh
  adaptation is essential for efficient computational solution of
  problems of the type discussed in \S\S \ref{sec:AL}--\ref{sec:ANISO}
  and, to the best of our knowledge, implementation of mesh adaptation
  in Ipopt-type, or more generally, in Newton-type methods, remains an
  open problem. We also emphasize that the (RCG) approach has far
  fewer parameters than the (Ipopt) method which greatly simplifies
  its performance optimization.} Finally, as a challenging test, the
(RCG) method was used to compute vortex configurations in rotating BEC
with high values of the nonlinear interaction constants and very high
rotation rates {as well as in configurations with strongly
  anisotropic trapping potentials.}

{Lastly, we reiterate that the approach presented in this study
  does not exploit all opportunities inherent in the Riemannian
  formulation. In particular, it remains an open question whether the
  use of a well-adapted Riemannian metric defined on the constraint
  manifold could further improve the performance of the approach. In
  addition, one can also consider the Riemannian formulations of
  Newton's method and of different variants of the quasi-Newton
  method.  Work is already on-going on some of these problems and
  results will be reported in the near future.}

\bigskip {The authors are grateful to the anonymous referees for
  providing constructive feedback on this manuscript.}

%


\end{document}